\documentclass[twocolumn,floats,amssymb,aps]{revtex4b1}
\usepackage{times,amsthm,amsmath,cite,pstricks,pst-node}

\makeatletter

\renewcommand{\p@subsection}{}

\renewcommand{\p@subsubsection}{}

\makeatother

\newtheorem{theorem}{Theorem}
\newtheorem{corollary}[theorem]{Corollary}
\newtheorem{lemma}[theorem]{Lemma}
\newtheorem{conjecture}[theorem]{Conjecture}
\theoremstyle{remark}
\newtheorem*{example}{Example}
\newtheorem*{remark}{Remark}
\theoremstyle{definition}

\renewcommand{\tensor}{\otimes}
\newcommand{\fin}{{\mathrm{fin}}}
\newcommand{\ie}{\emph{i.e.}}

\newcommand{\dis}{\mathrm{dis}}
\newcommand{\bub}{\mathrm{bub}}

\newcommand{\Pl}{\mathrm{Pl}}

\newcommand{\Th}{\mathrm{Th}}
\newcommand{\lk}{\mathrm{lk}}

\renewcommand{\d}{\partial}
\renewcommand{\bar}{\overline}
\renewcommand{\hat}{\widehat}
\newcommand{\sm}{\setminus}

\newcommand{\R}{{\mathbb{R}}}
\newcommand{\Z}{{\mathbb{Z}}}
\newcommand{\C}{{\mathbb{C}}}
\newcommand{\Q}{{\mathbb{Q}}}
\renewcommand{\H}{{\mathbb{H}}}
\newcommand{\bC}{\bar{C}}

\newcommand{\bP}{\bar{P}}
\newcommand{\bX}{\bar{X}}
\newcommand{\hC}{\hat{C}}
\newcommand{\hP}{\hat{P}}
\newcommand{\cC}{\mathcal{C}}
\newcommand{\cG}{\mathcal{G}}
\newcommand{\cS}{\mathcal{S}}
\newcommand{\cP}{\mathcal{P}}
\newcommand{\cM}{\mathcal{M}}
\newcommand{\mA}{\mathbf{A}}
\newcommand{\mB}{\mathbf{B}}
\newcommand{\mC}{\mathbf{C}}
\newcommand{\mD}{\mathbf{D}}
\newcommand{\mE}{\mathbf{E}}
\newcommand{\mF}{\mathbf{F}}
\newcommand{\mG}{\mathbf{G}}
\newcommand{\mH}{\mathbf{H}}
\newcommand{\mI}{\mathbf{I}}
\newcommand{\mX}{\mathbf{X}}
\newcommand{\tI}{\widetilde{I}}
\newcommand{\bul}{\bullet}
\newcommand{\del}{\Delta}

\newcommand{\eq}[2]{\begin{equation}\label{#1}#2\end{equation}}

\newcommand{\eatline}{\vspace{-\baselineskip}}
\newcommand{\fig}[1]{Figure~\ref{#1}}
\newenvironment{fullfigure}[2]
    {\begin{figure}[htb]\begin{center}\def\ffa{#1}\def\ffb{#2}}
    {\vspace{\baselineskip}\caption{\ffb.}\label{\ffa}\end{center}\end{figure}}

\newcommand{\mto}{\lput{:U}{\pspicture(0,0)(0,0)
    \psline[arrows=->,arrowscale=1.3](3.2pt,0)(3.3pt,0)\endpspicture}}
\newgray{gray15}{.85}
\psset{linewidth=.4pt,dash=5pt 5pt,unit=.25in}
\SpecialCoor

\begin{document}
\title{Perturbative 3-manifold invariants by cut-and-paste topology}
\author{Greg Kuperberg}
\thanks{Supported by NSF grant DMS \#9704125 and by a
    Sloan Foundation Research Fellowship}
\affiliation{UC Davis; {\tt greg@math.ucdavis.edu}}
\author{Dylan P. Thurston}
\thanks{Supported by an NSF Graduate Fellowship and a Sloan Foundation 
Dissertation Year Fellowship}
\affiliation{UC Berkeley; {\tt dpt@math.berkeley.edu}}
\begin{abstract}
  We give a purely topological definition of the perturbative quantum
  invariants of links and 3-manifolds associated with Chern-Simons
  field theory. Our definition is as close as possible to one given by
  Kontsevich. We will also establish some basic properties of these
  invariants, in particular that they are universally finite type with
  respect to algebraically split surgery and with respect to Torelli
  surgery.  Torelli surgery is a mutual generalization of blink
  surgery of Garoufalidis and Levine and clasper surgery of Habiro.
\end{abstract}
\maketitle

\section{Introduction}
\label{s:intro}

The purpose of this article is to give a purely topological definition of the
perturbative quantum invariants of links and 3-manifolds that were originally
defined by Kontsevich~\cite{Kontsevich:feynman} and that are associated with
Chern-Simons field theory. We will also establish some basic properties of
these invariants, in particular that they are finite type in the expected way.
The main difference between our treatment and that of Kontsevich is that we
will use cohomology rather than differential forms and the pairing between
homology and cohomology rather than integration. We will define the
perturbative invariants of degree $n$ of a given closed, framed rational
homology 3-sphere $M$ as the degree of a generalized Gauss map
$$\Phi:X_n \to Y_n,$$
where $X_n$ and $Y_n$ are modified configuration spaces that are constructed
from $M$ using cut-and-paste topology.  Otherwise we will follow Kontsevich
closely, since his definition is terse but essentially rigorous.

The purely topological approach was first considered by Bott and
Taubes~\cite{BT:knots}, and later by the second
author~\cite{Thurston:thesis} and
others~\cite{Poirier:rationality,AF:universal}.  The definition given
there generalized the Gauss map
$$\Phi: K_1 \times K_2 \to S^2,$$
whose degree is the linking number between two knots $K_1$ and $K_2$ in $\R^3$,
to other maps whose degrees give all of the Vassiliev invariants of knots and
links. Our maps can also be defined for links in 3-manifolds, a generalization
which we will discuss later.

More precisely, we will construct an invariant
\begin{align}
\Phi:(C_n,D) &\to (P^{\times 3n},Q) \\
\Phi^*:H^{6n}(P^{\times 3n},Q;\Q) &\to H^{6n}(C_n,D;\Q),
\label{e:universal}
\end{align}
where $6n$ is the degree of the top non-vanishing rational cohomology (or
homology) of certain spaces $C_n$ and $P^{\times 3n}$ that depend on $M$, and
$Q$ and $D$ are certain degenerate loci associated with the infinite asymptote
in an asymptotically flat model of $M$. The space $P$ has a generating class
$$\alpha \in H^2(P;\Q)$$
called a propagator.  The space $C_n$ is defined using the combinatorics of
Jacobi diagrams.
The space $V^*_n$ of primitive weight systems of degree $n$
embeds in the homology space $H_{6n}(C_n;\Q)$.  If $w \in V^*_n$
is a weight system, let $\mu_w$ be a corresponding cycle.
Then we can define an invariant
\eq{e:pairing}{I_w(M) = \langle w, \Phi^*(\alpha^{\tensor 3n}) \rangle 
= \langle \Phi_*(w), \alpha^{\tensor 3n} \rangle}
depending on a weight system.  Dually, we can define a universal
invariant (in a sense given below) as an element $I_n(M) \in V_n$.

\begin{theorem} The invariant $I_n(M)$ of framed rational homology
spheres is additive under connected sums:
\[I_n(M_1 \# M_2) = I_n(M_1) + I_n(M_2).\]
\eatline\label{th:sum}\end{theorem}

In particular, $I_n(S^3) = 0$ (if the modified tangent bundle $T'S$
as defined in Section~\ref{s:confblow} is given the canonical framing).

\begin{theorem} The invariant $I_n(M)$ is a finite-type invariant of degree $n$
in both the algebraically split and Torelli senses for framed rational homology
spheres $M$, and it is universal for integer homology spheres.
\label{th:fintype}
\end{theorem}

The phrase ``finite-type invariant'' merits some explanation.  In
general, suppose that $\cM$ is some set of topological objects with
the structure of a cubical complex $\cC$: Certain pairs of elements are
connected by edges, certain pairs of pairs form squares, and so on.
Then a function $I$ on $\cM$ (a topological invariant) taking values
in an abelian group extends to $\cC$ by taking alternating sums, or
repeated finite differences.  For example, if $M_{\emptyset}$, $M_1$,
$M_2$, and $M_{12}$ form a square $C$, then we can define
$$I''(C) = I(M_{\emptyset}) - I(M_1) - I(M_2) + I(M_{12}).$$
(We assume a suitable decoration on cubes to resolve the sign
ambiguity.) In this general context an invariant $I$ is
\emph{finite-type} of order $n$ if the $n+1$st order finite difference
$I^{(n+1)}$ vanishes. Another view is to interpret $S$ as a symbol for
a formal linear combination
$$S = M_{\emptyset} - M_1 - M_2 + M_{12}$$
and then extend $I$ linearly.  In this interpretation, we define $\cM_n$
to be the span (in the space of rational linear combinations)
of all $n$-cells of $\cC$.

As a motivating example, let $\cM_n$ be the set of $n$-dimensional
parallelepipeds on a vector space (i.e., a collection of $n$ vectors
together with a base point).  Then the functions on the vector space
that satisfy the above definition of finite-type of degree $n$ are
the polynomials of degree $n$.

In our case, $\cM$ is the set of homeomorphism types of oriented rational
homology 3-spheres, and the cubes in $\cC$ are defined in one of two ways:  The
vertices may be connected by surgery on sublinks of an algebraically split
link, or by Torelli surgery on subsets of a collection of disjoint
handlebodies.   Here an algebraically split link is a framed link whose linking
matrix is the identity; at the end of Section~\ref{s:fintype} we will also
consider a rational generalization.  A \emph{Torelli surgery} is the operation
of removing a handlebody from $M$ and gluing it back after applying an element
of the Torelli group to the boundary. (The Torelli group of a surface is the
subgroup of the mapping class group that acts trivially on homology.) 
Algebraically split surgery was defined by Ohtsuki~\cite{Ohtsuki:finite}, while
Torelli surgery generalizes both blink surgery as defined by Garoufalidis and
Levine~\cite{GL:blinks} and clasper surgery as defined by Habiro
\cite{Habiro:claspers}. Garoufalidis and Levine~\cite{GL:blinks} showed that
these two notions of finite type are equivalent to each other for integer
homology spheres. Moreover, they showed that there is a surjection
\eq{e:surject}{\kappa:V_n \twoheadrightarrow \cM_{kn}/\cM_{kn+1},}
where $k=3$ in the algebraically split case and $k=2$ in the Torelli case.  A
finite-type invariant is \emph{universal} if its finite difference of order
$kn$ is a right inverse
$$I^{(kn)}:\cM_{kn} \to V_n,$$
thereby showing that the map $\kappa$ is an isomorphism. We will argue
universality directly in both cases.  (Note that for unframed 3-manifolds,
$\cM_{kn+j} = \cM_{kn+j+1}$ when $k$ does not divide $j$.  The framed theory is
the same except that $\cM_1$ is 1-dimensional and detects change of framing.)

Finally in Section~\ref{s:unframed} we will prove the following theorem.

\begin{theorem} There is an invariant $\delta_n(M) \in V_n$ of 
homology 3-spheres $M$ decorated with a framing or bundled
bordism such that the difference
$$\tI_n(M) = I_n(M) - \delta_n(M)$$
is independent of the decoration.  Moreover, the framing
correction $\delta_n(M)$
is finite type of degree 1 in the Torelli and algebraically split
senses.
\label{th:correction}
\end{theorem}

In particular, the unframed invariant $\tI_n(M)$ is also universal.

It is also known~\cite{GGP:equivalence} that surgery on boundary links
(links whose components admit disjoint Seifert surfaces) again gives
an equivalent finite-type theory for integer homology spheres.  Thus
Theorem~\ref{th:fintype} has the following corollary.

\begin{corollary} The invariant $I_n(M)$ is universally finite type of order
$n$ for boundary link surgery in the class of integer homology spheres.
\label{c:boundary}
\end{corollary}

We do not yet have a direct proof of Corollary~\ref{c:boundary}.
We also have the following closely related conjectures.

\begin{conjecture} The unframed invariant
$$\omega(M) = \sum_n m^n \tI_n(M),$$
where $m = |H_1(M;\Z)|$, equals the surgery-defined invariant of Le, Murakami,
and Ohtsuki \cite{LMO:universal}. \label{c:lmo} \end{conjecture}

Conjecture~\ref{c:lmo} asserts that $\tI_n$ satisfies the Le-Murakami-Ohtsuki
surgery formula.  At the moment, we can only compute appropriate finite
differences to find the highest order term, analogous to the leading
coefficient of a polynomial.  Since both invariants are universal,
Conjecture~\ref{c:lmo} holds to highest order.

\begin{conjecture} The framing correction $\delta_n(M)$ vanishes
for $n>1$. \label{c:vanishes}
\end{conjecture}

By a remark in Section~\ref{s:remarks}, Conjecture~\ref{c:vanishes}
holds for $n$ even.

\subsection{Related work and further directions}

These definitions and results generalize to arbitrary rational homology spheres
and to links in rational homology spheres.  One interesting variant that we
have not analyzed is the definition of Axelrod and
Singer~\cite{AS:chern1,AS:chern2}, further developed by Bott and
Cattaneo~\cite{BC:integral1,BC:integral2}. The main difference between that
definition and the one due to Kontsevich (and ours) is that Kontsevich
punctures the 3-manifold $M$ so that the space of pairs of distinct points in
$M$ (the building block of the space $P$ above) is a homology 2-sphere, while
Axelrod and Singer ``smear out'' the puncture using a volume form. These
variations were considered in more detail by
Cattaneo~\cite{Cattaneo:configuration}. In this article we use a compact
version of Kontsevich's space, denoted $C_{e,\infty}(M)$; without puncturing it
would be just $C_e(M)$.  Algebraically, we need to know that
$$H^2(C_{e,\infty}(M);\Q) \cong \Q.$$
Following Axelrod and Singer, one could, without puncturing $M$, choose a
propagator
$$\alpha \in Z^2(C_e(M);\Q)$$
such that the coboundary satisfies
$$\delta\alpha = \mu \tensor 1 + 1 \tensor \mu,$$
where $\mu$ is a cocycle in $Z^3(M)$.

Taubes~\cite{Taubes:theta1,Taubes:theta2} defines and studies an invariant that
is very close to the invariant $I_1(M)$ that we define, using the canonical
framing of a rational homology 3-sphere.  He finds that his quantity is
invariant under spin cobordism, implying that it is trivial for integer
homology 3-spheres.  On the other hand, Theorems~\ref{th:fintype}
and~\ref{th:correction} imply that our invariant is the Casson invariant.  (A
standard relation for the  Casson invariant~\cite{AM:casson} implies that it is
finite type of degree 3 in the algebraically split sense; on the other hand the
space of invariants of this degree is 1-dimensional~\cite{Ohtsuki:finite}.)  We
have no explanation for the discrepancy, but we plan to consider the invariant
$I_1(M)$ in more detail in a future article~\cite{CHKT:copout}.

Two other generalizations that can be considered are invariants of graphs in
3-manifolds, and invariants associated to other flat
connections~\cite{BC:integral2}.  We will analyze these in future work.  Among
other things, there should be a general relation between flat bundles and links
in 3-manifolds on the one hand and finite covers and branched covers on the
other hand~\cite{Garoufalidis:cyclic}.

Kontsevich has discussed yet other generalizations.   There should be
corresponding invariants for a higher-dimensional smooth, framed manifold $M$
which produce certain characteristic classes of an $M$-bundle over another
topological space~\cite{Kontsevich:feynman}.  Our analysis may extend to
these invariants.  (Although the methods are still combinatorial, they do use
the tangent bundle, so it's not clear if the invariants would descend to PL
invariants.) More recently~\cite{Kontsevich:operads}, he explained that all
perturbative invariants are examples of homotopy functors from a certain
category of coordinate patches in $M$.

A more exotic possible generalization would be to pass from three real
dimensions to three complex dimensions.  It is possible that the holomorphic
cohomology of a Calabi-Yau 3-fold has all of the necessary properties to
generalize the definition of the invariant $I_n$.

\begin{acknowledgments}
We would like to especially thank Alberto Cattaneo and Nathan Habegger, who
helped us develop the theory and made several important comments.  We would
also like to thank Dror Bar-Natan, Andrew Casson, Rob Kirby, Paul Melvin,
and Bill Thurston for useful conversations.
\end{acknowledgments}

\section{Homological conventions}

None of the ideas in this article depend in any fundamental away on the model
of homology used: De Rham, singular, simplicial/cellular, \v{C}ech, etc.  For
concreteness it is convenient to use simplicial homology with coefficients in
$\Q$ and with unspecified triangulations.  We will henceforth omit the
coefficients.  (Note that most of the constructions would work identically with
arbitrary coefficients.)

Recall that the cup product in simplicial homology depends on an ordering of
the vertices, and that it is not graded commutative on the level of chains. 
This deficiency can be ameliorated when working over $\Q$ or any other
coefficient ring that contains $\Q$.  Namely, we can average over the $(a+b+1)!$
orderings of the vertices of each $a+b$-simplex when taking the cup product of
an $a$-cochain and a $b$-cochain.  If such a cochain is a cocycle, then it can
be represented by a differential form which is constant on each simplex, and a
cup product is then identically equal to the corresponding wedge product. In
this sense, simplicial cohomology is a kind of ``mock De Rham cohomology''.

The degenerate locus $Q$ is constructed as a semi-algebraic set rather than
with a cut-and-paste method.  Hence it does not a priori have simplicial
homology.  A foundational result of Hironaka states that semi-algebraic sets
can be ambiently triangulated, and the simplices of such a triangulation can be
straightened~\cite{Hardt:subanalytic,Hironaka:subanalytic}.

We will need the following extension lemma, which is elementary in the setting
of simplicial cohomology.

\begin{lemma} If $K$ is a subcomplex of a simplicial complex $L$, and if a
cohomology class $\alpha \in H^*(K)$ extends to a class $\beta \in H^*(L)$,
then any simplicial cocycle in $K$ representing $\alpha$ extends to a cocycle in $L$
representing $\beta$.
\label{l:extension}
\end{lemma}

\section{Jacobi diagrams}

In this section we review the definition of different kinds of Jacobi diagrams,
which are also variously called chord diagrams, Chinese characters, Chinese
character diagrams, and Feynman diagrams.  Technically we will need this
formalism only much later (in Lemma~\ref{l:domcohom} and
Section~\ref{s:properties}), but we present it here as a fundamental
preliminary.

\subsection{Parity functors}

Let $\cP$ be the category of two-element sets in which morphisms
are bijections; it has a natural tensor product operation if you view it as the
category of affine spaces over the multiplicative group $\{1,-1\}$.  (More
concretely:  The identity map from any object of $\cP$ to itself is called
$1$ and the other map is called $-1$. If $A = \{a,b\}$ and $X = \{x,y\}$ are in
$\cP$, then $A \tensor X$ has the two elements $\{(a,x),(b,y)\}$ and
$\{(b,x),(a,y)\}$.)  A parity functor is a functor from some other category with
invertible morphisms to $\cP$.   For example, let $\mA(S)$ be the set
of sign-orderings of a finite set $S$, \ie, the set of linear orderings
quotiented by the action of the alternating group $\mathrm{Alt}(S)$. $\mA$
and the orientation functor for finite-dimensional vector spaces are the two
most commonly used non-trivial parity functors.  (Arguably the trivial functor
$\mathbf{1}$, a special case of which is defined below, is even more commonly used.)

Let $\cG$ be the category of connected, finite graphs $\Gamma$ (multiple
edges and loops are allowed) in which the morphisms are graph isomorphisms.
One can consider the following parity functors on $\cG$:

\begin{description}
\item[$\mathbf{1}(\Gamma)$] is the trivial functor that takes every
graph to $\{1,-1\}$.
\item[$\mA(\Gamma)$] is the set of sign-orderings of the edges of $\Gamma$.
\item[$\mB(\Gamma)$] is the set of sign-orderings of the odd-valence vertices.
\item[$\mC(\Gamma)$] is the set of sign-orderings of the even-valence vertices.
\item[$\mD(\Gamma)$] is the set of orientations of all edges, quotiented by the
operation of negating any two.
\item[$\mE(\Gamma)$] is the set of sign-orderings of the flags (pairs
consisting of an edge and one of its vertices).
\item[$\mF(\Gamma)$] is the set of sign-orderings of the edges incident to each
vertex, up to negating any two sign-orderings.
\item[$\mG(\Gamma)$] is the set of sign-orderings of all vertices, equivalently
the set of orientations of the vector space of simplicial 0-chains
$Z_0(\Gamma;\R)$.
\item[$\mH(\Gamma)$] is the set of orientations of the vector space of 1-chains
$Z_1(\Gamma;\R)$.
\item[$\mI(\Gamma)$] is the set of orientations of $H_1(\Gamma;\R)$.
\end{description}

These functors, modulo isomorphism of functors (via natural transformations),
generate an abelian group with exponent 2 (since $\mX \tensor \mX \cong
\mathbf{1}$ for any parity functor $\mX$) with the relations:
\begin{align*}
\mE &\cong \mD             & \mH &\cong \mA \tensor \mE \\
\mF &\cong \mB \tensor \mE & \mI &\cong \mG \tensor \mH \\
\mG &\cong \mB \tensor \mC
\end{align*}
For example, the functors $\mD$ and $\mE$ are isomorphic as
follows:  An orientation of an edge $e$ of a graph $\Gamma$ can be expressed as an
ordering of the two flags that include $e$.  Listings the edges
$e_1,e_2,\ldots,e_n$ in any order, we get an ordering of the flags
$$f_{1,1},f_{1,2},f_{2,1},f_{2,2},f_{3,1},f_{3,2},\ldots,f_{n,2},$$
where the flags of the edge $e_i$ are ordered $(f_{i,1},f_{i,2})$. The sign of
this ordering of the flags does not depend on the ordering of the edges,
establishing a canonical isomorphism $\mD(\Gamma) \cong
\mE(\Gamma)$. We leave the other relations as an exercise.

Each parity functor determines a homomorphism
$$\mathrm{Aut}(\Gamma) \to \{1,-1\}.$$
There are choices for $\Gamma$ for which $\mA$, $\mB$, $\mC$, and
$\mD$ induce independent homomorphisms, for example the one in
\fig{f:independent}.

\begin{fullfigure}{f:independent}
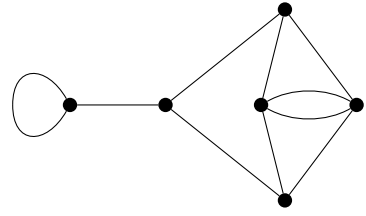
    {A graph demonstrating independence of parity functors
    $\mA$, $\mB$, $\mC$, and $\mD$}
\pspicture(-1,-2)(6,2)
\cnode*(0,0){.15}{f} \cnode*(2,0){.15}{e} \cnode*(4.5,-2){.15}{a}
\cnode*(4,0){.15}{c} \cnode*(6,0){.15}{b} \cnode*(4.5,2){.15}{d}
\pnode(-1.2,0){g}
\ncline{a}{b} \ncline{a}{c} \ncline{b}{d} \ncline{c}{d} \ncline{a}{e}
\ncline{d}{e} \ncline{e}{f}
\nccurve[angleA=30,angleB=150]{c}{b}
\nccurve[angleA=-30,angleB=-150]{c}{b}
\nccurve[angleA=120,angleB=90,ncurv=1.5]{f}{g}
\nccurve[angleA=270,angleB=240,ncurv=1.5]{g}{f}
\endpspicture
\end{fullfigure}

Thus, no further relations are possible.  The parity functors listed above can
be expressed in terms of the four generators according to Table~\ref{t:parity}.

\begin{table}
$$\begin{array}{c|ccccccccc}
     & \mA  & \mB  & \mC  & \mD  & \mE  & \mF  & \mG  & \mH  & \mI  \\ \hline
\mA  & \bul &      &      &      &      &      &      & \bul & \bul \\
\mB  &      & \bul &      &      &      & \bul & \bul &      & \bul \\
\mC  &      &      & \bul &      &      &      & \bul &      & \bul \\
\mD  &      &      &      & \bul & \bul & \bul &      & \bul & \bul
\end{array}$$
\caption{Nine parity functors in terms of four generators.}
\label{t:parity}
\end{table}

On the subcategory of $\cG$ of odd-valence graphs, $\mC \cong \mathbf{1}$, but $\mA$, $\mB$, and $\mD$ remain
independent.

We define a Lie orientation of $\Gamma$ to be an element of $(\mD \tensor
\mG)(\Gamma)$. This parity functor is naturally associated to
invariants and characteristic classes of odd-dimensional manifolds.  In the
association between graph homology and the twisted equivariant homology of
``outer space''~\cite{CV:moduli,Kontsevich:feynman}, the
isomorphic parity functor $\mA \tensor \mI$ appears. The parity functor $\mF
\tensor \mC$ is also isomorphic; Bar-Natan~\cite{Bar-Natan:vassiliev} defines
Lie orientations in the odd-valence case using just $\mF$.  Note
that the parity functor $\mA$
leads to the other kind of graph homology; it corresponds to the untwisted
equivariant homology of outer space and to configuration spaces
on even-dimensional manifolds.

\subsection{Diagrams and relations}

A \emph{closed Jacobi diagram} is a Lie-oriented graph $\Gamma$ with trivalent
vertices.  (A non-closed Jacobi diagram may also have univalent vertices.) A
closed diagram has $2n$ trivalent vertices if and only if it has $n+1$
\emph{loops}, where the loop number is just the first Betti number of the
diagram. The \emph{Vassiliev space} $V_n$ is the vector space over $\Q$ of
isomorphism classes of connected Jacobi diagrams with $n+1$ loops, modulo the
Jacobi relation (also called the $IHX$ relation):
\eq{e:Jacobi}{\pspicture[.4](-.85,-.7)(.85,.7)
\psline(-.7,.7)(0,.4)(.7,.7) \psline(-.7,-.7)(0,-.4)(.7,-.7)
\psline(0,-.4)(0,.4) \qdisk(0,.4){.08} \qdisk(0,-.4){.08} \endpspicture
= \pspicture[.4](-.85,-.7)(.85,.7)
\psline(.7,-.7)(.4,0)(.7,.7) \psline(-.7,-.7)(-.4,0)(-.7,.7)
\psline(-.4,0)(.4,0) \qdisk(.4,0){.08} \qdisk(-.4,0){.08} \endpspicture
- \pspicture[.4](-.85,-.7)(.85,.7)
\psline(-.7,-.7)(.7,.7) \psline(.7,-.7)(-.7,.7) \psline(-.35,-.35)(.35,-.35)
\qdisk(-.35,-.35){.08} \qdisk(.35,-.35){.08} \endpspicture}
This is a linear relation among any three graphs that are the same except at
the indicated subgraphs.  The edges incident to each vertex are cyclically
ordered clockwise in the diagram.  This is an $\mF$-orientation (the
\emph{blackboard} orientation), which by previous considerations is equivalent
to a Lie orientation for trivalent graphs.

We will also consider dual vectors $w \in V^*_n$, which are
called \emph{primitive weight systems}.

\begin{remark} The IHX relation is compatible with many kinds of decorations on
Jacobi diagrams.  The edges may be ordered; the homology or the fundamental
group may have distinguished elements or other decorations; there may be
univalent vertices which may or may not be labelled; and the diagram may be
attached to a link or a graph.  These decorations are important for
generalizations of the invariant $I_n(M)$ and for analyzing Vassiliev spaces,
but in this article we only need the simplest of all Vassiliev spaces.
\end{remark}

\begin{example} The spaces $V_n$ are 1-dimensional for $n=1,2$.
For $n=1$ there is only one diagram, the theta graph:
$$
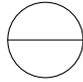
\pspicture(-1,-1)(1,1)\pscircle(0,0){.8}\psline(-.8,0)(.8,0)\endpspicture$$
For $n=2$ there are two, a double theta and a tetrahedron, and the
former is twice the latter:
$$
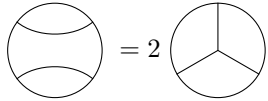
\pspicture[.43](-1.2,-1)(1.2,1)\pscircle(0,0){1}
\pcarc[arcangle=-40](-.8,.6)(.8,.6)
\pcarc[arcangle=40](-.8,-.6)(.8,-.6)\endpspicture
= 2\pspicture[.43](-1.2,-1)(1.2,1)\pscircle(0,0){1}
\psline(1;90)(0,0)(1;210)\psline(0,0)(1;330)\endpspicture
$$
As above, we assume the blackboard orientation in this equation.
\end{example}

\section{Configuration spaces}

In this section we will define a certain compactification of the configuration
space of maps from the vertices of a graph $\Gamma$ to a manifold $M$ such that
vertices connected by an edge are distinct. The idea is to blow up diagonals 
corresponding to the edges in the space of all maps $M^\Gamma$.  This is more
complicated than one might expect, since these diagonals are not mutually
transverse. We will rely on a general construction for resolving non-transverse
blowups of this type.

\subsection{Blowups:  The balls, beams, and plates construction}
\label{s:bbp}

In this section we will discuss blowing up a manifold $M$ along a general type
of closed subset $X$ called a \emph{Whitney-stratified space}
\cite{Thom:ensembles,GM:morse}.  By virtue of its Whitney stratification, $X$
decomposes into a locally finite, partially ordered set of smoothly embedded
manifolds,
$$X = \bigcup_{i \in \cS} X_i.$$
The decomposition and the partial ordering are compatible according
to the condition that 
$$i \prec j \iff X_i \subset \bar{X_j} \iff X_i \cap \bar{X_j} \ne \emptyset$$
for $i \ne j$.  In our case, we additionally require that $X$ is locally
smoothly equivalent to a cone over another Whitney-stratified space; \ie, for
each $p \in X$ there is a \emph{tangent cone} $T_pX$. We call such a
Whitney-stratified space \emph{cone-like}; one which is not conelike can have
cusps and other singularities in which strata kiss.

We will need a generalization of this definition which we call a
\emph{Whitney-stratified immersion}.  As before, $X$ decomposes into smoothly
embedded manifolds, and we assume that 
$$i \prec j \iff X_i \subset \bar{X_j}.$$
But the third condition, that $X_i$ and $\bar{X_j}$ are disjoint if $i$ and $j$
are incomparable, is replaced by two weaker conditions:
\begin{description}
\item[1.] Each $\bar{X_i}$ is a union of strata.
\item[2.] If $i_1,\ldots,i_n$ are an anti-chain,
then the corresponding strata $X_{i_1},\ldots,X_{i_n}$ are mutually transverse.
\end{description}
Here an \emph{anti-chain} is a pairwise incomparable set.

If $M$ is a manifold with a cone-like, Whitney-stratified immersion $X$, there
is a way to blow up $M$ along $X$.  It is convenient (but not strictly
necessary) to give $M$ a Riemannian metric. The blowup $B_X(M)$ is formed by
successively blowing up $X_i$ as $i$ increases. This means that we replace each
$p \in \bar{X_i}$ by the set of rays in $T_p(M)$ which are normal to
$T_p(X_i)$; here $\bar{X_i}$ is the closure of $X_i$ in the partially blown up
model of $M$.   If some strata in an anti-chain intersect transversely, then
their blowups commute, so they can be performed in either order.

The result $B_X(M)$ is a smooth manifold with right-angled corners: a manifold
locally diffeomorphic to a closed cube.  It has a codimension 1 face $F_i$ for
each $i$, and the interior of $F_i$ blows down to the open stratum $X_i$.
Lower-dimensional faces correspond to flags (ordered chains) of strata. This is
topologically and combinatorially equivalent to the complement of an open
regular neighborhood of $X$.  The latter is also called the ``balls, beams, and
plates'' construction when it appears in geometric topology.

\begin{example} Let $M$ be a square and let $X$ be a fish
on a line, as in \fig{f:fish}. 
Note that there would be a geometric pathology at the univalent points if we
tried to blow up all of $X$ in one go. These
pathologies become extreme in high dimensions.
\end{example}

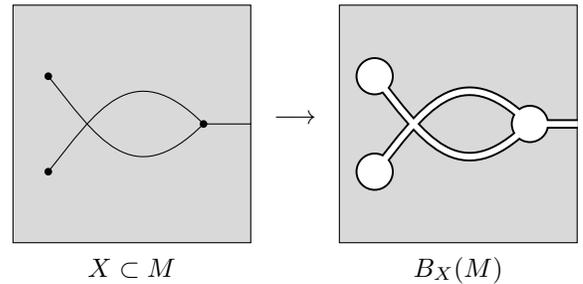
\begin{fullfigure}{f:fish}{Iterated blowup of a square at a fish on a line}
\begin{tabular}{ccc}
\pspicture[.5](-2.75,-2.75)(2.75,2.75)
\psframe[fillstyle=solid,fillcolor=gray15](-2.5,-2.5)(2.5,2.5)
\psbezier(1.5,0)(0,1.5)(-.75,.25)(-1.75,-1)
\psbezier(1.5,0)(0,-1.5)(-.75,-.25)(-1.75,1)
\psline(1.5,0)(2.5,0)
\qdisk(1.5,0){.08}
\qdisk(-1.75,-1){.08}
\qdisk(-1.75,1){.08}
\endpspicture & $\longrightarrow$ &
\pspicture[.5](-2.75,-2.75)(2.75,2.75)
\psframe[fillstyle=solid,fillcolor=gray15](-2.5,-2.5)(2.5,2.5)
\psbezier[linewidth=.2](1.5,0)(0,1.5)(-.75,.25)(-1.75,-1)
\psbezier[linewidth=.2](1.5,0)(0,-1.5)(-.75,-.25)(-1.75,1)
\psline[linewidth=.2](1.5,0)(2.5,0)
\qdisk(1.5,0){.4}
\qdisk(-1.75,-1){.4}
\qdisk(-1.75,1){.4}
\psbezier[linewidth=.12,linecolor=white](1.5,0)(0,1.5)(-.75,.25)(-1.75,-1)
\psbezier[linewidth=.12,linecolor=white](1.5,0)(0,-1.5)(-.75,-.25)(-1.75,1)
\psline[linewidth=.12,linecolor=white](1.5,0)(2.5,0)
\pscircle*[linecolor=white](1.5,0){.36}
\pscircle*[linecolor=white](-1.75,-1){.36}
\pscircle*[linecolor=white](-1.75,1){.36}
\endpspicture \\
$X \subset M$ & & $B_X(M)$
\end{tabular}
\end{fullfigure}

\begin{remark} The construction is actually more general in several important
respects.  First, instead of a Whitney-stratified immersion in a manifold, we
could consider an immersion of one cone-like space in another one.  Even if the
target space is a manifold, this  allows the blowup locus $X$ to be a
transversely immersed submanifold, for example.  Second, in our blowups we
quotient by multiplication by scalars in $\R^+$.  We could instead quotient by
scalar multiplication by $\R^\ast$ or $\C^\ast$, provided that the tangent cone
at each point in $X$ is invariant under this larger group of homotheties.
The iterated $\C^\ast$ blowup of complex configuration spaces is called the 
Fulton-Macpherson compactification~\cite{FM:compactification}.
\end{remark}

\subsection{Geometry of blowups}

If $p \in F_k$ blows down to $q$ for $q \in X_i$,
then $p$ can be thought of as a point ``infinitely close''
to $q$.  More formally, it is an element of the quotient $(T_q(M) -
T_q(X_i))/\R_+$, where $\R_+$ acts by positive rescaling in the
directions normal to $T_q(X_i)$.  We can
and will use the vector space structure of $T_q(M)$ to describe $p$.  If $p$
lies in a corner of $B_X(M)$, for example in the intersection of $F_i$
and $F_j$ for $i < j$, then it can be understood as infinitely close to both
$X_i$ and $X_j$, but infinitely closer to $X_j$ than
to $X_i$.

In the main construction we will label part of the blowup locus as being ``at
infinity'' and give $T_p(M)$ an inverted linear structure at points $p$ in this
locus. In the simplest example, $M$ has a marked point $\infty$.  Define
$M_\fin$, the finite part of $M$, as
$$M_\fin = M \sm \{\infty\}.$$
If $M$ has a Riemannian metric, then we can give $M_\fin$ an asymptotically
flat Riemannian metric by inverting the exponential map from the point
$\infty$. If we add a sphere at infinity to $M_\fin$ in the usual way by adding
endpoints to infinite rays, the result $\bar{M_\fin}$ is combinatorially
equivalent to the blowup $B_\infty(M)$ of $M$ at $\{\infty\}$.

Although it is more complicated to describe, in the general case any subcomplex
$Y \subset X$ can be considered the infinite locus, and the blowup of $M$ along
$Y$ can be given an inverted geometry.  The idea is to invert the exponential
map normal to each stratum $Y_i$.

These geometries will be described more explicitly in the case of interest in
the next section.

\subsection{Blowups for configuration spaces}
\label{s:confblow}

Let $M$ be a $d$-dimensional manifold and let $\Gamma$ be a connected graph
with $n$ vertices.  The graph $\Gamma$ may have self-loops and multiple edges,
but these do not affect the construction in this section.  Let the symbol
$\Gamma$ also denote the vertex set of $\Gamma$, so that
$$M^\Gamma = \{f:\Gamma \to M \}$$
is equivalent to a Cartesian product $M^{\times n}$.  Points in $M^{\Gamma}$,
and in other spaces that we will form from it, are called
\emph{configurations}. If $\Gamma'$ is a subgraph of $\Gamma$, let
$\Delta_{\Gamma'}$ denote the diagonal in $M^\Gamma$ in which all vertices of
$\Gamma'$ are sent to the same point.

In order to define a Gauss map, we need to blow up the diagonal $\Delta_e$ for
every edge $e \subset \Gamma$.  Such a diagonal is called \emph{principal}, and
blowing it up produces a codimension 1 face $F_e$ called a \emph{principal
face}. The principal diagonals are not mutually transverse, so we must blow up
other diagonals first.  Specifically, we blow up $\Delta_{\Gamma'}$ for every
vertex-2-connected subgraph $\Gamma' \subseteq \Gamma$.  (A \emph{cut vertex}
of a connected graph is a vertex whose removal disconnects the graph.  A graph
is \emph{vertex-2-connected} if it is connected and has no cut vertices.  Note
that a single edge is vertex-2-connected.)  These will be used as the strata in
Section~\ref{s:bbp}. Unless $\Gamma'$ is a single edge, the diagonal
$\Delta_{\Gamma'}$ is called \emph{hidden} and blowing it up produces a
codimension 1 face $F_{\Gamma'}$ called a \emph{hidden face}.  If $\Gamma$
itself is 2-connected, then $F_{\Gamma}$ is called the \emph{anomalous face}. 
We denote the result $C_\Gamma(M)$, the (compactified) $\Gamma$-configuration
space of $M$. Its construction is valid modulo the following lemma:

\begin{lemma} The system of diagonals $\{\Delta_{\Gamma'}\}$ corresponding to
2-connected subgraphs forms a conelike, Whitney-stratified, self-transverse
immersion in $M^\Gamma$.
\label{l:transverse}
\end{lemma}
\begin{proof} (Sketch) Checking that the diagonals are Whitney-stratified and
conelike is complicated but routine; the more significant issue is
self-transversality.

We describe the minimal strata containing a configuration $c \in M^{\Gamma}$.
Let $\Gamma_1,\ldots,\Gamma_k$ be the connected components of the inverse
images in $c$ of points in $M$, not including components which consist of a
single point. If $c$ lies on some of the diagonals, then this list of subgraphs
is non-empty.  A subgraph $\Gamma_i$ might not be 2-connected. Rather, it is
has tree-like structure consisting of maximal 2-connected subgraphs which share
cut vertices.  We call such a structure a \emph{cactus} and the 2-connected
subgraphs \emph{lobes}. A diagonal $\Delta_{\Gamma'}$ is a minimal stratum
containing $c$ if and only if $\Gamma'$ is a maximal 2-connected subgraph of
some $\Gamma_i$.  Checking that these strata are mutually transverse is again
complicated but routine.
\end{proof}

Each codimension 1 face $F_{\Gamma'}$ has a geometric structure that we will
use to describe certain gluings. The face $F_{\Gamma'}$ fibers over a smaller
configuration space $C_{\Gamma/\Gamma'}(M)$, where $\Gamma/\Gamma'$ is the
graph $\Gamma$ with $\Gamma'$ contracted to a vertex $p \in M$. Let
$f_{\Gamma',p}$ be a fiber where no other point of $\Gamma$ is close to $p$. 
Then this fiber is just $C_{\Gamma'}(T_pM)/\Th$, where $\Th = \Th(T_pM)$ is the
$d+1$-dimensional Lie group of translation and homothety (scalar multiplication
by $\R_+$) in the tangent space $T_pM$. Later we will need the quotient
$$c_{\Gamma}(V) = C_{\Gamma}(V)/\Th(V)$$
for an arbitrary $d$-dimensional vector space $V$. If $\Gamma' = e$ is an edge,
each $f_{e,p}$ is diffeomorphic to $S^{d-1}$. The reader can check that 
$C_\Gamma(M)$ is $dn$-dimensional and that $dn-1$ is the total dimension
of the fibration structure of each codimension 1 face.

The general corner (\ie, face with codimension $\ge 2$) of
$C_\Gamma(M)$ is given by a list of codimension 1 faces that meet it,
or equivalently a list of edges and 2-connected subgraphs.  In general
corners in the balls, beams, and plates construction can either come
from transverse intersections or from flags of strata; corners of both
types are illustrated in Figure~\ref{f:fish}.  In our case corners
that are purely of the first type have the same combinatorics as the
corresponding transverse intersections as described in the proof of
Lemma~\ref{l:transverse}:  As a configuration approaches the corner, the
graph $\Gamma$ develops one or more cactus structures whose lobes are
2-connected subgraphs; the vertices in each node converge and the nodes
converge together. Corners that are purely of the second type consist
of nested 2-connected subgraphs
$$\Gamma_1 \supset \Gamma_2 \supset \ldots \supset \Gamma_k,$$
where possibly the innermost graph $\Gamma_k$ is a single edge. In this
case the vertices of each $\Gamma_i$ draw together as a configuration
approaches the corner, but those of $\Gamma_{i+1}$ draw together  at
a faster rate than those of $\Gamma_i$.  Most corners are of mixed
type:  Each is given by a cactus forest of 2-connected subgraphs,
and there may be more cactus forests nested in the cactus lobes.

\begin{example} Figure~\ref{f:codim4} shows an example of a configuration $c$
in which nine vertices of $\Gamma$ have converged.  Because these vertices form
a double square and a triangle connected by an edge, $c$ lies on a principal
face and two hidden faces.  In addition two of the vertices in the double
square have converged more quickly to each other than to the other four
vertices in their cluster, which means that $c$ lies on a second principal face
as well.  Thus, the face has codimension 4.  (Note that some of the gluings
described in Section~\ref{s:gluings} reduce the dimensions of some of the faces
and corners.)
\end{example}

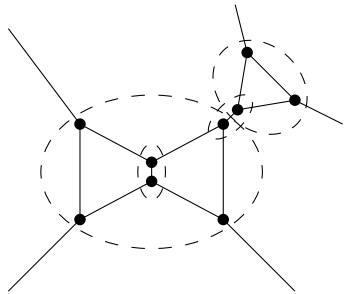
\begin{fullfigure}{f:codim4}{A nested cactus configuration on a codimension 4 face}
\pspicture(-3.5,-2.5)(4,3.5)
\cnode*(-1.5, 1){.12}{a}
\cnode*( 1.5, 1){.12}{b}
\cnode*(-1.5,-1){.12}{c}
\cnode*( 1.5,-1){.12}{d}
\cnode*( 0, .2){.12}{e}
\cnode*( 0,-.2){.12}{f}
\cnode*(1.8,1.3){.12}{g}
\cnode*(2,2.5){.12}{h}
\cnode*(3,1.5){.12}{i}
\psellipse[linestyle=dashed](0,0)(2.33,1.62)
\psellipse[linestyle=dashed](0,0)(.3,.6)
\rput{-45}(1.65,1.15){\psellipse[linestyle=dashed](0,0)(.3,.6)}
\rput{45}(2.27,1.77){\psellipse[linestyle=dashed](0,0)(.88,1.1)}
\psline(e)(b)(d)(f)(c)(a)(e)(f)
\psline(g)(h)(i)(g)(b)
\psline(a)(-3,3)
\psline(c)(-3,-2.5)
\psline(d)( 3,-2.5)
\psline(h)(1.75,3.5)
\psline(i)(4,1)
\endpspicture
\end{fullfigure}

In the main construction of the article, we work with $M$ a closed manifold
with a marked point $\infty \in M$.   We will need a compact configuration
space $C_{\Gamma,\infty}(M)$ which contains the partially compactified space
$C_{\Gamma,\infty}(M_\fin)$. For each 1-connected subgraph $\Gamma' \subseteq
\Gamma$, we blow up $M^\Gamma$ along the locus $\Delta_{\infty,\Gamma'}$ where
all vertices of $\Gamma'$ lie at the point $\infty$.  (It is not enough to do
this only for 2-connected subgraphs $\Gamma'$.  In effect $\Gamma$ has been
suspended from a special vertex $\infty$, and $\infty \cup \Gamma'$ is
2-connected if and only if $\Gamma'$ is 1-connected.) We treat this locus as an
infinite locus in the sense of the previous section.  These blowups together
with the diagonal blowups yield the space $C_{\Gamma,\infty}(M)$.

\begin{example}
Suppose for simplicity that we only blew up $M^\Gamma$ for one such subgraph
$\Gamma'$ and that we did not blow up $\Delta_{\Gamma''}$ for any subgraph
$\Gamma''$ of $\Gamma'$.   Let $R^{\Gamma'}$ be the set of endpoints of rays
in $M_{\fin}^{\Gamma'}$, which is asymptotically flat just as $M_{\fin}$ is.
Then the blown up locus is $R^{\Gamma'} \times M_{\fin}^{\Gamma - \Gamma'}$.
In a configuration in the blown up locus, the vertices of $\Gamma'$ lie at an
astronomical scale compared to $M_{\fin}$.  We retain the relative distances
between these points and $M_{\fin}$ and the angles, but not the scale of these
distances relative to the internal geometry of $M_{\fin}$. See \fig{f:mfin}.
\end{example}

Finally, to match the inverted geometry at $\infty$, we will use a modified
tangent bundle $T'M$.  This is the unique bundle on $M$ whose sections pull
back to asymptotically constant sections on $TM_\fin$; alternatively, it is the
push forward of the bundle $T(M \# B^3)$ to $M$, mapping the boundary sphere of
$B^3$ to a point and identifying the fibers of $TB^3$ using its natural
trivialization (from the inclusion $B^3 \subset \R$). If $M$ is an orientable
3-manifold, $T'M$ is isomorphic to $TM$ since they are both trivial, but the
isomorphism is not canonical.  (Note that $T'S^d$ is always trivial, while
$TS^d$ is non-trivial for even $d$.)

\begin{fullfigure}{f:mfin}{A configuration at the astronomical scale}
\pspicture(-.5,-.5)(5,5)
\pspolygon(0,0)(5,2)(2,5)
\qdisk(5,2){.12} \qdisk(2,5){.12}
\qdisk(3,3){.12}
\psline(5,2)(3,3)(2,5)
\psline(3,3)(0,0)
\cput[fillstyle=solid,fillcolor=white](0,0){$\scriptstyle M_{\fin}$}
\endpspicture
\end{fullfigure}

\section{Existence of the invariant}

\subsection{The gluings}
\label{s:gluings}

In this section we will construct a grand configuration space $C_n$ using all
decorated, closed Jacobi diagrams $\Gamma$ with $n+1$ loops. Multiple edges are
allowed, but self-loops are not. The way that we decorate the diagrams is that
we explicitly orient every edge and we fully order the edges.  We will need the
Lie orientation of $\Gamma$ only later, in Lemma~\ref{l:domcohom}.  The
vertices of $\Gamma$ are not ordered. Let $\hat{J}_n$ be the set of such
diagrams up to isomorphism.

Let $M$ be a closed 3-manifold with a marked point $\infty$, and assume a
framing of the modified tangent bundle $T'M$.  This is equivalent to an
asymptotically constant framing of $M_\fin$ and is not much different from a
framing of $M$.

We start with a dismembered version of the Gauss map $\Phi$. Let $C_n(M)$ be the
disjoint union of
all $C_{\Gamma,\infty}(M)$:
$$C_n(M) = \coprod_{\Gamma \in \hat{J}_n} C_{\Gamma,\infty}(M)$$
and let
$$P(M) = C_{e,\infty}(M)$$
where the graph $e$ is an edge.
The map
$$\Phi:C_n(M) \to P(M)^{\times 3n}$$
is defined in the $k$th factor by erasing the vertices of $\Gamma$ other than
the two in the $k$th edge.  Thus the configuration space $P(M)$ is a kind of
topological propagator.  It has the desired homology $H^2(P(M)) \cong \Q$, but
$C_n(M)$ has no degree $6n$ (or top) homology because each component has
faces.  So we will glue the faces of $C_n(M)$ to each other, or otherwise cap,
collapse, or relativize them, and correspondingly modify $P(M)$ as necessary
without destroying its homology.  This process is equivalent to Kontsevich's
arguments that certain improper integrals vanish or cancel. The result will be
a commutative diagram:
$$\begin{array}{c@{\hspace{1cm}}c}
\rnode{a}{C_n(M)} & \rnode{b}{P(M)^{\times 3n}} \\[1cm]
\rnode{c}{\bC_n(M)} & \rnode{d}{\bP(M)^{\times 3n}}
\end{array}
\ncline[nodesep=.3]{->}{a}{b}\Aput{\Phi}
\ncline[nodesep=.3]{->}{c}{d}\Aput{\Phi}
\ncline[nodesep=.25]{->>}{a}{c}
\ncline[nodesep=.25]{->>}{b}{d}
$$
As desired, the Vassiliev space $V_n$ will appear as a quotient of the top
cohomology of the glued configuration space $\bC_n(M)$ computed
relative to the the degenerate locus $D$.

The gluings are as follows:

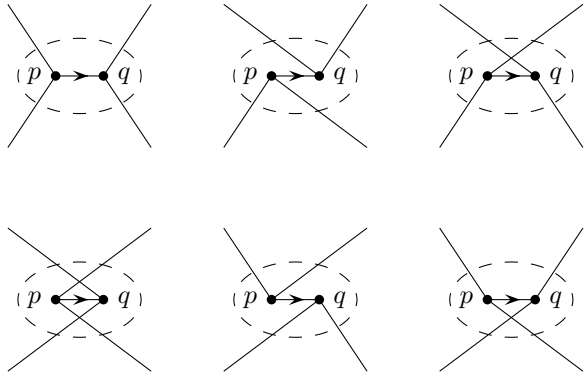
\begin{fullfigure}{f:six}{Identifying six principal faces}
\begin{tabular}{ccc}
\pspicture(-2.15,-2.25)(2.15,1.5)
\psline(-1.5,-1.5)(-.5,0)(-1.5,1.5) \psline(1.5,-1.5)(.5,0)(1.5,1.5)
\pcline(-.5,0)(.5,0)\mto            \qdisk(-.5,0){.1} \qdisk(.5,0){.1}
\rput[r](-.8,0){$p$} \rput[l](.8,0){$q$}
\psellipse[linestyle=dashed](0,0)(1.3,.8)
\endpspicture &
\pspicture(-2.15,-2.25)(2.15,1.5)
\psline(-1.5,-1.5)(-.5,0)(1.5,-1.5) \psline(-1.5,1.5)(.5,0)(1.5,1.5)
\pcline(-.5,0)(.5,0)\mto            \qdisk(-.5,0){.1} \qdisk(.5,0){.1}
\rput[r](-.8,0){$p$} \rput[l](.8,0){$q$}
\psellipse[linestyle=dashed](0,0)(1.3,.8)
\endpspicture &
\pspicture(-2.15,-2.25)(2.15,1.5)
\psline(-1.5,-1.5)(-.5,0)(1.5,1.5) \psline(1.5,-1.5)(.5,0)(-1.5,1.5)
\pcline(-.5,0)(.5,0)\mto           \qdisk(-.5,0){.1} \qdisk(.5,0){.1}
\rput[r](-.8,0){$p$} \rput[l](.8,0){$q$}
\psellipse[linestyle=dashed](0,0)(1.3,.8)
\endpspicture \\
\pspicture(-2.15,-1.5)(2.15,2.25)
\psline(-1.5,-1.5)(.5,0)(-1.5,1.5) \psline(1.5,-1.5)(-.5,0)(1.5,1.5)
\pcline(-.5,0)(.5,0)\mto           \qdisk(-.5,0){.1} \qdisk(.5,0){.1}
\rput[r](-.8,0){$p$} \rput[l](.8,0){$q$}
\psellipse[linestyle=dashed](0,0)(1.3,.8)
\endpspicture &
\pspicture(-2.15,-1.5)(2.15,2.25)
\psline(-1.5,-1.5)(.5,0)(1.5,-1.5) \psline(-1.5,1.5)(-.5,0)(1.5,1.5)
\pcline(-.5,0)(.5,0)\mto           \qdisk(-.5,0){.1} \qdisk(.5,0){.1}
\rput[r](-.8,0){$p$} \rput[l](.8,0){$q$}
\psellipse[linestyle=dashed](0,0)(1.3,.8)
\endpspicture &
\pspicture(-2.15,-1.5)(2.15,2.25)
\psline(-1.5,-1.5)(.5,0)(1.5,1.5)  \psline(1.5,-1.5)(-.5,0)(-1.5,1.5)
\pcline(-.5,0)(.5,0)\mto           \qdisk(-.5,0){.1} \qdisk(.5,0){.1}
\rput[r](-.8,0){$p$} \rput[l](.8,0){$q$}
\psellipse[linestyle=dashed](0,0)(1.3,.8)
\endpspicture \\
\end{tabular}
\end{fullfigure}

\begin{description}
\item[Principal faces:] For each simple edge $e$ we glue together the principal
faces $F_{e}$ of the configuration spaces $C_{\Gamma,\infty}(M)$ for six
different graphs $\Gamma$.   The graphs are paired by reversing the
orientation of $e$, and the three pairs differ by the Jacobi relation with
$e$ in the middle, as in \fig{f:six}.  The ordering of the edges
changes only implicitly, by virtue of the fact that the edges are
reconnected.  (This rule implies that six distinct faces are always
glued.)

\item[Hidden faces:] Recall that the hidden faces correspond to 2-connected
subgraphs $\Gamma'$; a double edge is counted as a hidden face rather than a
principal one. Suppose that a pair of edges $e_1$ and $e_2$ in $\Gamma'$
separates it into two subgraphs $\Psi_1$ and $\Psi_2$.  Then we can glue
$F_{\Gamma'}$ to another face in which  $e_1$ and $e_2$ are switched and their
orientations are reversed.  Since all points of  $\Gamma'$ lie in some tangent
space $T_pM$, modulo the homothety group, we can describe this operation
explicitly in terms of linear algebra in $T_pM$. We leave
$\Psi_2$ fixed and send every vertex $q \in \Psi_1$ to $q-e_1-e_2$, where
$e_1$ and $e_2$ also denote vectors corresponding to the edges $e_1$ and
$e_2$ point from $\Psi_2$ to $\Psi_1$.  Note that this involution changes
the extra decoration on $\Gamma$, namely the ordering and orientation of its
edges, but not the underlying graph.

We need to know that there is at least one involution. since $\Gamma'$ is not
all of $\Gamma$, it has a vertex $q$ connected to $\Gamma\setminus\Gamma'$. 
Since $\Gamma'$ is 2-connected, $q$ has valence 2 in $\Gamma'$.  We let
$\Psi_1$ be $q$ and call the neighboring vertices $p_1$ and $p_2$, so that $e_1
= (p_1,q)$ and $e_2 = (p_2,q)$.  (If $\Gamma'$ is a double edge, then $p_1 =
p_2$.) \fig{f:hidden} illustrates the hidden face involution in this case.

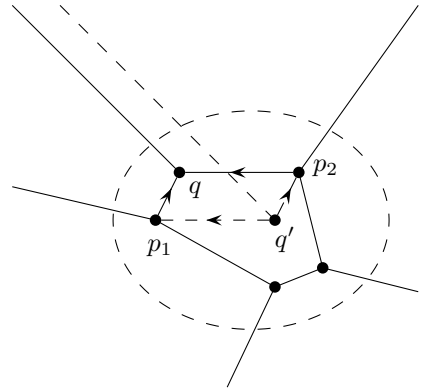
\begin{fullfigure}{f:hidden}{The involution for a hidden face}
\pspicture(-4.5,-4)(4,4)
\cnode*(-1.5,-.5){.12}{a}   \cnode*(-1,.5){.12}{b}
\cnode*(1.5,.5){.12}{c}     \cnode*(1,-.5){.12}{d}
\cnode*(1,-1.9){.12}{e}
\cnode*(2,-1.5){.12}{f}
\rput[t](-1.4,-.9){$p_1$}  \rput[tl](-.8,.3){$q$}
\rput[l](1.8,.6){$p_2$}    \rput[tl](1,-.65){$q'$}
\ncline{a}{b}\mto
\ncline{c}{b}\mto
\ncline[linestyle=dashed]{d}{a}\mto
\ncline[linestyle=dashed]{d}{c}\mto
\psline[linestyle=dashed](d)(-3.5,4)
\psline(b)(-4.5,4) 
\psline(-4.5,.2)(a)(e)(f)(c)(4,4)
\psline(0,-4)(e)
\psline(4,-2)(f)
\psellipse[linestyle=dashed](.5,-.5)(2.9,2.3)
\endpspicture
\end{fullfigure}

\item[The anomalous face:] This face is a compactification of a bundle with
fiber $c_{\Gamma}(T(M_\infty))$ over $\bar{M_\fin}$. We identify all
fibers with each other using the framing of $M_\fin$. We perform the same
operation in the topological propagator $P(M)$.

The unique face corresponding to $\Theta$, the unique Jacobi diagram with two
vertices, is treated as anomalous rather than principal or hidden.

\item[Infinite faces:]  First, the topological propagator $P(M)$ has two
semi-infinite faces with one vertex at infinity and the other not, and it has a
totally infinite face with both vertices at infinity.  A configuration  in any
of these faces determines an element of $S^2$ by taking the unit vector point
from vertex 1 to vertex 2.  We identify all three faces with standard $S^2$
using this correspondence; this $S^2$ is necessarily identified with the
remnant of the anomalous face of $P(M)$.  Denote the result $\bP(M)$.

The infinite faces of the domain $C_n(M)$, including the totally infinite
face, form the degenerate locus $D$.  The degenerate locus
$$Q = \bigcup_A Q_A \times {\bP(M)}^{\times 3n-|A|}$$
is the union of pieces, one for each set $A$ of the edges numbered from 1 to
$3n$. Given $A = \{a_1,\ldots,a_k\}$, the locus $Q_A$ consists of those
elements
$$(v_{a_1},\ldots,v_{a_k}) \in (S^2)^{\times A}$$
such that the unit vectors $\{v_{a_i}\}$ can be realized as the directions of
the edges of some graph with $k$ edges which has been linearly mapped into
$\R^3$. The graph is required to have no vertices of valence 1 and at most
one of valence 2, although multiple edges are allowed.

\end{description}

The result is a glued configuration space $\bC_n(M)$ and a glued
topological propagator $\bP(M)$.

\subsubsection{The transitive closure of the gluings}

To describe the transitive closure of the gluings of the
principal and hidden faces, we begin with the geometry of the configuration
space after only the principal faces are glued.  The principal edges of a
configuration $c$ form a forest of trees in the graph $\Gamma$.  (Because of
the hidden blowups, these edges cannot form closed loops.)  By virtue of the
principal blowups, which may be performed simultaneously after all hidden
blowups, each of these edges has a well-defined direction but not a length, not
even a relative length when compared with  any other edge of $\Gamma$.  The
principal gluings then  identify $c$ with all other configurations in which
each tree of principal edges is replaced by some other tree with the edges
pointing in the same directions.  In the glued space the trees lose their
identity.  The data that remains is a graph $\bar{\Gamma}$ in which each tree
of principal edges in $\Gamma$ is contracted to a point; a vertex in
$\bar{\Gamma}$ at the point $p \in M$ of valence $n>3$ is also assigned a list
of $n-3$ unit tangent vectors in $T_pM$.

Some of the corners that are glued to each other do not have the same
dimension, because the reconnection in \fig{f:six} can changed whether or not a
subgraph is 2-connected.  Section~\ref{s:confblow} describes how a corner
before gluing is determined by nested cactus structures in the diagram
$\Gamma$, and that the codimension of the corner equals the total number of
lobes of the cacti.  After the gluings of the principal faces, a general corner
is described by the contracted graph $\bar{\Gamma}$ together with nested cacti
in $\bar{\Gamma}$.  Each lobe of each cactus now corresponds to a hidden face
and cannot be a single edge.  The total codimension of the corner is then the
total number of principal edges in $\Gamma$ plus the total number of cactus
lobes in  $\bar{\Gamma}$.  For example, before the principal gluings the
configuration on the left in \fig{f:codim4} lies on a codimension 4 face.  The
reconnection in \fig{f:six} glues it to the codimension 5 face in
\fig{f:codim5}.

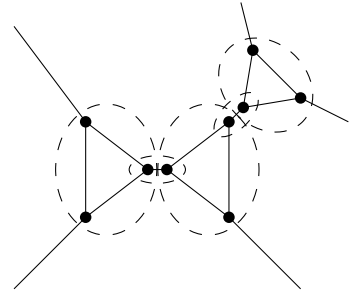
\begin{fullfigure}{f:codim5}{A cactus configuration on a codimension 5 face}
\pspicture(-3.5,-2.5)(4,3.5)
\cnode*(-1.5, 1){.12}{a}
\cnode*( 1.5, 1){.12}{b}
\cnode*(-1.5,-1){.12}{c}
\cnode*( 1.5,-1){.12}{d}
\cnode*(-.2,0){.12}{e}
\cnode*(.2,0){.12}{f}
\cnode*(1.8,1.3){.12}{g}
\cnode*(2,2.5){.12}{h}
\cnode*(3,1.5){.12}{i}
\psellipse[linestyle=dashed](0,0)(.6,.3)
\rput{-45}(1.65,1.15){\psellipse[linestyle=dashed](0,0)(.3,.6)}
\rput{45}(2.27,1.77){\psellipse[linestyle=dashed](0,0)(.88,1.1)}
\psellipse[linestyle=dashed](-1.07,0)(1.07,1.385)
\psellipse[linestyle=dashed](1.07,0)(1.07,1.385)
\psline(e)(a)(c)(e)(f)(b)(d)(f)
\psline(g)(h)(i)(g)(b)
\psline(a)(-3,3)
\psline(c)(-3,-2.5)
\psline(d)( 3,-2.5)
\psline(h)(1.75,3.5)
\psline(i)(4,1)
\endpspicture
\end{fullfigure}

Assuming that the principal gluings are completed, we describe the  hidden face
involutions.  The description rests on two facts.  First, following the
description in the previous paragraph, a hidden face involution always glues
together faces of the same dimension.  If two edges $e_1$ and $e_2$ of a
configuration of some face separate the contracted subgraph $\Gamma'$, then the
inverse image of a cactus lobe in $\bar{\Gamma}$ contains either both or
neither.  This is true even if one or both of $e_1$ or $e_2$ is principal.
Second, the involutions for any given $\Gamma$ generate a finite group,
because they are given by permuting and reversing edges.  (In fact
it is a product of symmetric groups, each acting on an equivalence class
of edges, where two edges are equivalent if they separate $\Gamma'$.)

\subsubsection{Remarks on the construction}
\label{s:remarks}

The entire construction has a folded version in which $P(M)$ is defined as the
space of unordered pairs of points, the Cartesian product $\bP(M)^{\times 3n}$
is replaced by the symmetric power $S^{3n}(\bP(M))$, and $\R_+$-blowups are
replaced by $\R^*$-blowups throughout.  The propagator space $P(M)$ becomes a
homology $\R P^2$ rather than a homology sphere, and its relevant second
cohomology group has coefficients in the twisted flat line bundle over $P(M)$
or $\bP(M)$.  The edges of $\Gamma$ are no longer explicitly oriented, nor are
the edges ordered.  This version is formally cleaner, but it is harder to
visualize.  It essentially hides signs and denominators in homological algebra
rather than removing them.

Many hidden faces admit more symmetries than those  generated by the given
involutions.  Namely for each subgraph $\Gamma'' \subset \Gamma'$ connected to
$\Gamma'$ at two vertices $p_1$ and $p_2$, we can reverse every edge in
$\Gamma''$ and switch $p_1$ and $p_2$ relative to $\Gamma'' \setminus \Gamma'$.
Also we can reverse every edge of $\Gamma'$. These involutions generate a
larger gluing group.  The key property of any such gluing group is that half of
its elements negate the map $f$ in Lemma~\ref{l:domcohom}.  Reversing
$\Gamma''$ has this effect if and only if $\Gamma''$ has an odd number of edges
plus vertices; reversing all of $\Gamma'$ does if and only if it has an even
number of edges plus vertices.  Thus reversing a single edge does not negate
$f$, but if $n$ is even, reversing all of $\Gamma$ in the anomalous face does
\cite[\S6]{AS:chern2}.  This removes the need to collapse the anomalous face
using the framing; in physics terminology, the anomaly cancels.  The other
involutions available also eliminate the anomaly if $\Gamma$ isn't
edge-3-connected.

Kontsevich \cite{Kontsevich:feynman} blows up every diagonal of $M^{\times 2n}$
to form an analogue of $C_\Gamma(M)$ that does not depend on $\Gamma$. While
his convention may have a rigorous analytic interpretation, it does not work
well in topologically, because a hidden face involution can send some vertices
of $\Gamma'$ on top of others.  In other works
\cite{Poirier:rationality,Thurston:thesis,BC:integral1} (and in the first
version of this work), the diagonal $\Delta_{\Gamma'}$ is blown up for every
connected subgraph $\Gamma' \subseteq \Gamma$.  Then the hidden face
involutions are more complicated and also inconsistent at corners.  The usual
remedy is to observe that corners are irrelevant to degree calculations; only
codimension 1 faces are important.  This is equivalent to relegating all of the
corners and their images to the relative locus.  We feel that it is more
natural to only blow up up $\Delta_{\Gamma'}$ when $\Gamma'$ is 2-connected.

Another approach to relativizing the semi-infinite faces, which
may be what Kontsevich had in mind, is power counting. If $\alpha \in
H^2(\bP(M))$ (defined in Section~\ref{s:cohomology}) is a Hodge form, it
vanishes as $L^{-2}$ on a length scale $L$ in the asymptotic part of
$\bP(M)$.  At the same time the available volume for a single vertex grows as
$L^3$.  The product is a negative power of $L$ for semi-infinite faces of
$C_\Gamma(M)$, which means that these faces are irrelevant in the degree
formula for $\Phi$.  It may be possible to phrase this argument in terms of
spectral sequences of filtrations, since $\bP(M)^{\times 3n}$ can be filtered
according to how many coordinates lie in $S^2 \subset \bP(M)$, while
$\bC_n(M)$ can be filtered according to how many vertices are at infinity.

The power counting argument does not work for the totally infinite face.  In
this case an alternative is to cap $F_{\Gamma,\infty}(M)$ with
$C_{\Gamma}(S^3)$, since the geometry of the face does not depend on the
manifold $M$.

\subsection{The bordism variant of framings}
\label{s:bordism}

Instead of collapsing the anomalous face of $C_\Gamma(M)$, we can instead cap
it using a bordism of $M$.  Although a special case of this formulation is
entirely equivalent to the framing approach, it will be more convenient for
the constructions in Section~\ref{s:fintype}.

More precisely, let $W$ be a 4-manifold bounded by $M$ and let $E$ be a
3-plane bundle that restricts to $T'M$ on $M$. Then we can cap the anomalous
face $F_\Gamma(M)$ with a certain configuration bundle $c_\Gamma(E)$ over $W$
for all graphs $\Gamma$.  The fiber over $p \in W$ of this bundle is the
configuration space $c_\Gamma(E_p)$ of the fiber $E_p$. This configuration
bundle has its own principal and hidden faces, which are glued in the same way
as faces of $C_{\Gamma,\infty}(M)$ to form a bundle $c_n(E)$. All of the other faces of
$C_{\Gamma,\infty}(M)$ are also glued the same way as before.  We denote the
resulting glued space $\hC_n(M)$.

Likewise we can cap the diagonal face of $P(M)$ with $c_e(E)$, which is just
the unit sphere bundle $SE$.   We can also refrain from collapsing the
semi-infinite faces or the totally face of $P(M)$. Call the result $\hP(M)$. 
This propagator space may have some spurious second cohomology coming from the
homology of $W$, but there is a unique second cohomology class in $SE$ which
can be represented by a cocycle which is antisymmetric under the antipodal map
on fibers.  (The antipodal map on the fibers, which extends to the map
switching the two factors of $C_2(M)$, splits the (rational)
cohomology into the odd and even subspace.  All the cohomology classes
from $W$ are even, by definition.)
This class extends to the propagator class $\alpha \in
\hP(M)$. The Gauss map $\Phi$ is defined as before.

If $T'M$ has a framing that extends to $E$, then there is a quotient map
$$\pi:\hC_n(M) \to \bC_n(M)$$
given by collapsing $W$ to a single point (and $E$ to a single fiber).  The map
$\pi$ induces an isomorphism of the top homology of the configuration spaces,
and the analogous map on propagators takes $\alpha$ to $\alpha$.  (This
property can be used as the definition of $\alpha \in \hP(M)$.)  The map $\pi$
then forms a commutative square with $\Phi$:
$$\begin{array}{c@{\hspace{1cm}}c}
\rnode{a}{\hC_n(M)} & \rnode{b}{\hP(M)^{\times 3n}} \\[1cm]
\rnode{c}{\bC_n(M)} & \rnode{d}{\bP(M)^{\times 3n}}
\end{array}
\ncline[nodesep=.3]{->}{a}{b}\Aput{\Phi}
\ncline[nodesep=.3]{->}{c}{d}\Aput{\Phi}
\ncline[nodesep=.25]{->>}{a}{c}\Bput{\pi}
\ncline[nodesep=.25]{->>}{b}{d}\Bput{\pi}
$$
This square and the isomorphism properties of $\pi$ demonstrate that, if the
bundle $E$ matches the framing of $T'M$, $\hC_n(M)$ produces the same
invariant $I_n(M)$ as $\bC_n(M)$.

Indeed, the bordism $W$ need not be a manifold, but only a homology manifold. 
(A homology $n$-manifold for us is a simplicial complex such that the link of
each vertex is a homology $n-1$-manifold and a homology $n-1$-sphere.)  In
particular, if $W$ is the cone over $M$, a bundle over $W$ extending $T'M$ is
equivalent to a framing of $T'M$.

\subsection{Cohomology}
\label{s:cohomology}

\begin{lemma} \label{l:domcohom}
 The top cohomology of the glued configuration space
$\bC_n(M)$ is independent of $M$ and has a surjection onto the
Vassiliev space $V_n$:
$$H^{6n}(\bC_n(M), D) \twoheadrightarrow V_n.$$
\end{lemma}

\begin{proof} Let $X_n$ be the union of all faces of $C_n$ (both finite
and infinite) and let $\bX_n$ be its image in $\bC_n$.   Consider the
cohomology exact sequence of the triple $D \subset \bX_n \subset \bC_n$:
$$H^{6n-1}(\bX_n,D) \to H^{6n}(\bC_n,\bX_n) \to H^{6n}(\bC_n, D) \to 0.$$
On the other hand,
$$
H^{6n}(\bC_n,D\cup\bX_n) \cong H^{6n}(C_n, X_n) \cong
  \bigoplus_{\Gamma \in \hat{J}_n}\Q\Gamma
$$
since $D\cup\bX_n$ cuts $\bC_n$ into the pieces of $C_n$, and on each
of these we have a unique top cohomology class, the fundamental class.
Thus, $H^{6n}(\bC_n(M),D)$ is a space of graphs modulo some relations.
These graphs are not quite Lie-oriented as in the definition
of $V_n$, since the edges are labelled and each edge is oriented.
(A Lie-oriented graph has a global choice of orientations up to sign).
But there is a forgetful map $f$ from $\hat{J}_n$ to Lie-oriented graphs.

To prove the lemma, we only need to check that the relations given by
$H^{6n-1}(\bX_n, D)$ become trivial or the Jacobi relation under $f$.
The space $H^{6n-1}(\bX_n,D)$ might be rather
complicated, but by the same exact sequence it is generated
by one cohomology class for each face (of $C_n$).  By cases:
\begin{description}
\item[Principal Faces:] The sum of the six graphs in \fig{f:six} descends to 
the Jacobi relation.
\item[Hidden faces:] We glue together several different graphs which become
identical (up to sign) under $f$.  Each involution defined in
Sections~\ref{s:gluings} negates $f(\Gamma)$. For example, if the graph
$\Psi_2$ is the single point $q$, then
$$q \mapsto  p_1 + p_2 - q$$
is orientation-reversing, and two edges are reversed.  Thus
half of the elements of the group generated by these involutions
negates $f(\Gamma)$, so the total sum vanishes.
\item[Infinite faces:] Since cohomology is computed relative to
these faces, they impose no relation.
\item[Anomalous face:] Since we reduce the
dimension of this face, it imposes no relation.
\end{description}
\end{proof}

\begin{remark} It may appear as if we are discarding information present in the
rest of $H^{6n}(\bC_n(M), D)$ by relying on Lemma~\ref{l:domcohom}. However,
the arguments of Section~\ref{s:fintype} imply that any invariant determined by
the action of $\Phi$ on $H^{6n}(\bC_n(M), D)$ is finite type.  Since $I_n$ is
universal among finite type invariants by Theorem~\ref{th:fintype}, it
determines all other such invariants.  In the minimal construction mentioned in
Section~\ref{s:remarks}, $H^{6n}$ is isomorphic to $V_n$; the spurious
cohomology is absent.
\end{remark}

\begin{lemma}\label{l:propcohom}
If $M$ is a rational homology sphere, then
the second cohomology $H^2(\bP(M))$ of the glued topological
propagator is generated by the fundamental class $\alpha$ of the standard
sphere $S^2 \subset \bP(M)$. Moreover, there is a well-defined cohomology class
$$\alpha^{\tensor 3n} \in H^{6n}(\bP(M)^{\times 3n},Q).$$
\end{lemma}
\begin{proof} The existence of $\alpha$ originates with the geometry of the
configuration space $C_{e,\infty}(M)$. This is a manifold with corners whose
interior is $M_{\fin}^{\times 2} \sm \Delta$, the space of pairs of  distinct
points in $M_{\fin}$.  If $M_{\fin} = \R^3$, it is clearly homotopy equivalent
to $S^2$.  In the general case it has the same homology by a Mayer-Vietoris
argument. Each of the gluings used to make $\bP(M)$ from $P(M)$ is chosen
to preserve the second cohomology, although higher cohomology may also appear.

The class $\alpha^{\tensor 3n}$ clearly exists in the absolute cohomology
$H^{6n}(\bP(M)^{\times 3n})$; the question is whether it exists uniquely in
cohomology relative to $Q$.  Observe first that if $|A|=k$, then $Q_A \subset
(S^2)^k$ has codimension at least 3. Each allowed graph $\Gamma'$ on $A$ with
$k$ edges (of which there are finitely many) has at most $(2k+1)/3$ vertices.
Thus there are at most $2k+1$ degrees of freedom in embedding $\Gamma'$ in
$\R^3$. In addition, 4 of these degrees of freedom are absorbed by invariance
under the homothety group $\Th(\R^3)$, so $Q_A$ has dimension at most $2k-3$.

Choose a point
$$p = (p_1,\ldots,p_{3n}) \in (S^2)^{\times 3n} \sm Q.$$
For each $i$, choose a cocycle $\alpha_i \in Z^2(\bP(M))$ that represents the
class $\alpha$ and that is localized at $p_i$ (or for concreteness, a small
simplex containing $p_i$) in the standard sphere $S^2 \subset \bP(M)$.\ Recall
that the space of relative cocycles $Z^{6n}(\bP(M)^{\times 3n},Q)$ is a
subspace of the space of absolute cocycles $Z^{6n}(\bP(M))$. The cocycle
$$\alpha_p = \alpha_1 \tensor \ldots \tensor \alpha_{3n}$$
exists as a relative cocycle because it avoids $Q$.  It represents a
non-trivial cohomology class because relativization can only diminish the
space of boundaries.  Thus $\alpha^{\tensor 3n}$ exists in relative
cohomology.

To show uniqueness, suppose that $a \subset S^2$ is an arc
connecting $p_1$ with some point $p'_1$ and which is disjoint from $Q$:
$$a \times (p_2,p_3,\ldots,p_{3n}) \subset (S^2)^{\times 3n} \sm Q.$$
If $\alpha'_1$ represents $\alpha$ and is localized at $p'_1$, then there is a
1-cochain $\beta$ localized along $a$ which is a homology between $\alpha_1$
and $\alpha'_1$:
$$\delta \beta = \alpha'_1 - \alpha_1.$$
In this case
$$\beta \tensor \alpha_2 \tensor \alpha_3 \tensor \ldots \tensor \alpha_{3n}$$
is a homology between $\alpha_p$ and $\alpha_{p'}$, where
$$p' = (p'_1,p_2,p_3,\ldots,p_{3n}).$$
Since $Q$ has codimension 3 in $(S^2)^{\times 3n}$, any two points in its
complement can be connected by a sequence of moves of this type.  Hence
$\alpha^{\tensor 3n}$ is unique.
\end{proof}

Having defined all elements of the map~\eqref{e:universal} and
equation~\eqref{e:pairing} (taking $C_n = \bC_n(M)$ and $P = \bP(M)$), the
definition of the invariant $I_n(M)$ for framed, rational homology spheres is
complete.

\section{Properties}
\label{s:properties}

\subsection{Connected sums}

In this subsection we prove Theorem~\ref{th:sum}. Although the conclusion  is
in the spirit of properties of surgery, the argument has more in common with
the definition of $I_w(M)$.

Suppose that $M = M_1 \# M_2$ is a rational homology sphere.  We may realize
$M_\fin$ by patching very small copies of $(M_1)_\fin$ and $(M_2)_\fin$ into a
flat $\R^3$, as in \fig{f:sum}.  In fact, $(M_1)_\fin$ and $(M_2)_\fin$ can be
infinitely small.  More precisely, we blow up $\R^3$ (with the trivial framing)
at $(0,0,0)$ and $(1,0,0)$ and we glue the spheres at infinity of
$(M_1)_\infty$ and $(M_2)_\infty$ to the blown up points.  In addition to the
astronomical scale used to compactify $M_\fin$, this geometry gives it an
intermediate scale, which we called the \emph{planetary scale}. On the
planetary scale, $(M_1)_\fin$ and $(M_2)_\fin$ are reduced to points but
are a unit distance from each other.  Call the set of points in this region
$\Pl(M_1,M_2)$.

We use the planetary scale to compactify $P(M)$ slightly differently.  If
$(p,q) \in P(M)$ and at least one of $p$ and $q$ is in $\Pl(M_1,M_2)$, or if
one is in $(M_1)_\fin$ and the other is in $(M_2)_\fin$, we glue $(p,q)$ to
the point in $S^2$ given by the direction from $p$ to $q$.

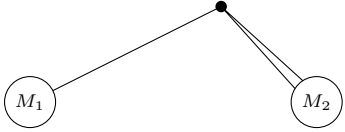
\begin{fullfigure}{f:sum}{The planetary scale for $M_1 \# M_2$}
\pspicture(-3.5,-.5)(3.5,2.1)
\qdisk(1,2){.12} \psline(1,2)(-3,0) \psline(3.2,0)(1,2)(2.8,0)
\cput[fillstyle=solid,fillcolor=white](-3,0){$\scriptstyle M_1$}
\cput[fillstyle=solid,fillcolor=white](3,0){$\scriptstyle M_2$}
\endpspicture
\end{fullfigure}

We may slightly enlarge the degenerate locus $Q$ without changing
$H^{6n}(\bP(M)^{\times 3n},Q)$.  In the definition of $Q$ in
Section~\ref{s:gluings}, we allow graphs in $\R^2$ with at most two vertices
of valence 2 rather than at most one, and we allow the graph consisting of a
single edge from $(0,0,0)$ to $(1,0,0)$ (or vice versa).  The resulting locus
$Q'$ has codimension 2 rather than codimension 3, but 2 is still
enough for the arguments of Lemmas~\ref{l:domcohom} and~\ref{l:propcohom}

If a configuration in $\bC_n(M)$ has any points at the planetary scale as in
\fig{f:sum}, or if it has some points in $(M_1)_\fin$ and others in
$(M_2)_\fin$, then the map $\Phi$ sends it to the locus $Q'$.  The only other
possibilities are that all vertices are in $(M_1)_\fin$, or that all vertices
are in $(M_2)_\fin$.  This realizes the cocycle $\alpha^{\tensor 3n}$
as the sum of cocycles on $\bC_n(M_1)$ and $\bC_n(M_2)$, which establishes
the identity
$$I_n(M) = I_n(M_1)+ I_n(M_2).$$

\subsection{Surgeries}
\label{s:surgeries}

In order to argue Theorem~\ref{th:fintype}, we would like to add and subtract
the cohomological propagators for different 3-manifolds.  Since these
propagators are defined on configuration spaces for different manifolds, we
will dismember the configuration spaces so that some of the pieces are the
same, and then calculate with the propagators on these common pieces.

We begin by more precisely defining the cubical complex $\cC$
mentioned in Section~\ref{s:intro} in the algebraically split
and Torelli cases.

We will consider a knot $K$ in a 3-manifold $M$ to be a closed solid torus that
does not contain the marked point $\infty$, and a link $L$ to be the union of
finitely many disjoint knots $\{K_1,\ldots,K_k\}$. For each such link $L$ we
will consider the $2^k$ sublinks of the form
$$L_I = \bigcup_{i \in I} K_i$$
where $I \subseteq [k] = \{1,\ldots,k\}$ is a set of indices. For each such
$L_I$ we will let $M_I$ be the result of $+1$ surgery on each component of
$L_I$.  Recall that a link $L$ in an integer homology sphere is
\emph{algebraically split} if the linking number between each pair of
components vanishes.  In the case we interpret the pair $(M,L)$ as an element
of $\cM$ given by the alternating sum
$$(M,L) = \sum_{I \subseteq [k]} (-1)^{|I|} M_I.$$
Thus if $f(M)$ is an invariant, then $f^{(k)}(M,L)$, the $k$th algebraically
split finite difference of $f$, is defined by the same sum.

We use the same conventions for Torelli surgery. As mentioned in the
introduction, a \emph{Torelli surgery} on an integer homology 3-sphere $M$
consists of removing a handlebody $H$ and gluing back a handlebody $H'$ that
differs by a surface automorphism which acts trivially on $H^1(\d H)$ (an
element of the Torelli group of $\d H$).  The locus $T$ of a Torelli surgery is
the union of finitely many disjoint handlebodies $\{H_1,\ldots,H_k\}$ (a
\emph{multi-handlebody} in $M$), where each $H_i$ is decorated with an element
of the Torelli group of $\d H_i$. For each multi-handlebody $T$ we will
consider the sub-multi-handlebodies $T_I$ for each $I \subseteq [k]$
and we let  $M_I$ be the result of surgery on $T_I$. We let 
$$(M,U) = \sum_{I \subseteq [k]} (-1)^{I} M_I$$
and if $f(M)$ is an invariant, we let $f^{(k)}(M,T)$ be the $k$th Torelli
finite difference of $f$.

Consider a surgery (either algebraically split or Torelli) on a manifold $M$ in
which a submanifold $N$ is replaced by some other submanifold $N'$. If $M$ is
framed, we will assume that $N'$ has a framing which agrees with the framing of
$N$ at the boundary. Likewise if $M$ has a bundle bordism $(W,E)$, then we
will assume a cobordism $W'$ between $N$ and $N'$ to attach to $W$. If $W$ has
a bundle $E$ extending the modified tangent bundle $T'M$, we can extend it to
$W'$.  The choices for this extra data will not matter, as long as we always
make the same choice for a surgery component $N$ which is shared by many
multi-component surgeries.  Note that because of spin obstructions,
a framing does not extend across an algebraically split surgery
if any of the surgery slopes are odd.  But bundle bordisms
always extend.

\subsection{Dismemberment and bubble wrap}
\label{s:gory}

The best way to understand dismemberment of a manifold $M$ is as a kind of
blowing up.  If $S$ is a surface in $M$, we can blow up $M$ along $S$, which
amounts to cutting $M$ along $S$, to make a manifold $B$.   We can also add
configurations in $M^\Gamma$ that meet $S$ to the blowup loci used to construct
$C_{\Gamma,\infty}(M)$. Call the resulting configuration space
$C_{\Gamma,\infty}(B)$. There is a blow-down map
$$\pi:C_{\Gamma,\infty}(B) \to C_{\Gamma,\infty}(M).$$
For example, suppose that $M$ consists of two manifolds $M_1$ and $M_2$
identified along a connected surface $S$.  (It is immaterial here which of
$M_1$ and $M_2$ has the point $\infty$, as long as it is not on $S$ itself.)
Then the blowup $Z$ is
$$B = M_1 \amalg M_2.$$
If $e$ is an edge, then  $C_e(B)$ has four components, defined by which of the
vertices of the edge are in $M_1$ and which are in $M_2$.  The four components
are homeomorphic to $C_e(M_1)$, $C_e(M_2)$, $M_1 \times M_2$, and $M_2 \times
M_1$.   Their geometry is slightly different, because if and $p,q \in M$ are
coincident on $S$, the point $(p,q)$ is blown up to record the direction from
$p$ to $q$ and the ratio of the distance from $p$ to $S$ to the distance from
$q$ to $S$.  Nonetheless by abuse of notation we will refer to the components
as $C_e(M_1)$, $M_1 \times M_2$, etc.

In the definition of $\bC_n(M)$, the gluings of the hidden faces and the
anomalous face are difficult to reconcile with blowing up along a surface $S$.
However, the anomalous face poses no problem if we cap it using a bundle
bordism $(W,E)$, since we can then extend $S$ to a hypersurface $T$ in $W$
and blow that up too.  Thus the topological propagator $\hP(M)$ can be
dismembered to make $\hP(B)$. Instead of dismembering $\hC_n(M)$, we will pull
back propagators defined on it to the pieces $C_{\Gamma,\infty}(M)$ and
$c_\Gamma(E)$, which we will then dismember.

In comparing propagators, we only need to compare the first
algebraically split discrete derivative.  Let $M$ be an integer
homology sphere and let $K = K_1 \subset M$ be a knot. Let $M_1$ be
the result of replacing $K$ by $K'$ in $M$, where $K'$ and $K$ differ
by a $+1$ Dehn twist.

\begin{lemma} If two integer homology spheres $M$ and $M_1$ differ by $+1$
surgery on a knot $K$, and if $\alpha \in H^2(\hP(M))$ and $\alpha_1 \in
H^2(\hP(M_1))$ are cohomological propagators, then $\alpha_1 - \alpha$ is
homologous to $\beta_1 \tensor \beta_1$ on $\hP(M\sm K) = \hP(M_K \sm K')$,
where $\beta_1$ is a 1-cocycle dual to a Seifert surface of $K$ (a Seifert
cocycle).
\label{l:difference}
\end{lemma}

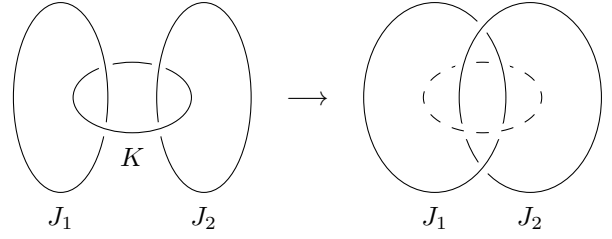
\begin{fullfigure}{f:link}{Linking two knots with unknot surgery}
\begin{tabular}{ccc}
\pspicture[.5](-3,-2.5)(3,2.2)
\psellipse(0,0)(1.25,.75)
\psellipse[linewidth=.3,linecolor=white](-1.5,0)(1.15,2.15)
\psellipse(-1.5,0)(1,2)
\psellipse[linewidth=.3,linecolor=white](1.5,0)(1.15,2.15)
\psellipse(1.5,0)(1,2)
\psclip{\psframe[linestyle=none](-3,-2.5)(3,0)}
\psellipse[linewidth=.3,linecolor=white](0,0)(1.4,.9)
\psellipse(0,0)(1.25,.75)
\endpsclip
\rput[t](-1.5,-2.3){$J_1$}
\rput[t](1.5,-2.3){$J_2$}
\rput[t](0,-1.05){$K$}
\endpspicture & $\longrightarrow$ & 
\pspicture[.5](-3,-2.5)(3,2.2)
\psellipse[linestyle=dashed](0,0)(1.25,.75)
\psellipse[linewidth=.3,linecolor=white](-1,0)(1.65,2.15)
\psellipse(-1,0)(1.5,2)
\psellipse[linewidth=.3,linecolor=white](1,0)(1.65,2.15)
\psellipse(1,0)(1.5,2)
\psclip{\psframe[linestyle=none](-3,-2.5)(3,0)}
\psellipse[linewidth=.3,linecolor=white](-1,0)(1.65,2.15)
\psellipse(-1,0)(1.5,2)
\psellipse[linewidth=.3,linecolor=white](0,0)(1.4,.9)
\psellipse[linestyle=dashed](0,0)(1.25,.75)
\endpsclip
\rput[t](-1,-2.3){$J_1$}
\rput[t](1,-2.3){$J_2$}
\endpspicture
\end{tabular}
\end{fullfigure}

\begin{proof} We can measure $\alpha_1 -\alpha$ by pairing it with 2-cycles in
$\hP(M\sm K)$.  There are several kinds of these, but the only kind that can
have non-zero pairing is represented by a torus $J_1 \times J_2$, where $J_1$
and $J_2$ are two disjoint knots in $M\sm K$. In this case $\alpha$ measures
their linking number in $M$:
$$\langle J_1 \times J_2 , \alpha \rangle = \lk_M(J_1,J_2)$$
Likewise $\alpha'$ measures their linking number in $M_K$. The
difference is the product of linking numbers with $K$:
$$\lk_{M_1}(J_1,J_2) - \lk_M(J_1,J_2) = \lk_M(J_1,K)\lk_M(J_2,K)$$
This is easy to see when $K$ is an unknot, since surgery on $K$ has the effect
of twisting $J_1$ and $J_2$ about each other without changing $M$, as in
\fig{f:link}. Since $\alpha_1-\alpha$ pairs with homology classes in the same
way as $\beta_1 \tensor \beta_1$, the two cocycles are homologous.
\end{proof}

The significance of Lemma~\ref{l:difference} is that by
Lemma~\ref{l:extension}, we can define $\alpha_1$ to be an extension of
$\alpha$ adjusted by $\beta_1$:
\eq{e:beta}{\alpha_1 \stackrel{\mathrm{def}}{=}
\alpha + \beta_1 \tensor \beta_1}
on $\hP(M\sm K)$.
Note also that we can assume that the support of $\beta_1$ is a neighborhood
of any desired Seifert surface $S$ of $K$.

The next case is algebraically split surgery with two components.  Let $L =
\{K_1,K_2\}$ be a link in $M$. Then the each of the four topological
propagators $\hP(M)$, $\hP(M_1)$, $\hP(M_2)$, and $\hP(M_{1,2})$ dismember into
nine pieces. The dismemberment of $\hP(M)$ looks like this:
$$\pspicture(-5.25,-2)(5.25,2)
\rput(-3.5,1.2){$\hP(K_1)$}\rput(0,1.2){$K_1\times M\sm L$}
\rput(3.5,1.2){$K_1\times K_2$}
\rput(-3.5,0){$M\sm L \times K_1$}\rput(0,0){$\hP(M\sm L)$}
\rput(3.5,0){$M\sm L \times K_2$}
\rput(-3.5,-1.2){$K_2 \times K_1$}\rput(0,-1.2){$K_2 \times M\sm L$}
\rput(3.5,-1.2){$\hP(K_2)$}
\psframe[linestyle=dashed,framearc=.5](-5.25,-.6)(1.75,1.8)
\psframe[linestyle=dashed,framearc=.5](-1.75,.6)(5.25,-1.8)
\endpspicture$$
Here we have circled $\hP(M \sm K_1)$ and $\hP(M \sm K_2)$.
Choosing Seifert surfaces $S_1$ and $S_2$ and Seifert cocycles $\beta_1$ and
$\beta_2$, we define $\alpha_1$ and $\alpha_2$ by equation~\eqref{e:beta} and
the extension principle. We assume that $S_1$ is disjoint from $K_2$ and vice
versa. Finally $\hP(M_{K_1,K_2})$ dismembers as follows:
$$\pspicture(-5.25,-2)(5.25,2)
\rput(-3.5,1.2){$\hP(K'_1)$}\rput(0,1.2){$K'_1\times M\sm L$}
\rput(3.5,1.2){$K'_1\times K'_2$}
\rput(-3.5,0){$M\sm L \times K'_1$}\rput(0,0){$\hP(M\sm L)$}
\rput(3.5,0){$M\sm L \times K'_2$}
\rput(-3.5,-1.2){$K'_2 \times K'_1$}\rput(0,-1.2){$K'_2 \times M\sm L$}
\rput(3.5,-1.2){$\hP(K'_2)$}
\psframe[linestyle=dashed,framearc=.5](-5.25,-.6)(1.75,1.8)
\psframe[linestyle=dashed,framearc=.5](-1.75,.6)(5.25,-1.8)
\endpspicture$$
In this diagram the northwest square is shared with $\hP(M_1)$, while the
southeast square is shared with $\hP(M_2)$.  By the boundary-disjointness of
the Seifert surfaces, if we define
$$\alpha_{K_1,K_2} = \alpha +
\beta_1 \tensor \beta_1 + \beta_2 \tensor \beta_2$$
on $\hP(M\sm L)$, we can extend it by $\alpha_1$ and $\alpha_2$ on the
rest of the shared pieces.  This leaves the two remaining pieces $K'_1 \times
K'_2$ and $K'_2 \times K'_1$. We claim that $\alpha_{1,2}$ automatically
extends to these pieces, because they can cannot create any second homology. In
other words, the inclusion
$$\hP(M_{1,2}) \sm ( K'_1 \times K'_2 \cup K'_2 \times K'_1) \subset \hP(M_{1,2})$$
is an isomorphism on $H^2$. This may be seen by a general position argument,
where we abbreviate the inclusion as just $X \subset Y$: Since $K'_1 \times
K'_2$ and $K'_2 \times K'_1$ are thickened 2-tori in the interior of $Y$, a
6-manifold with boundary, any 2-cycle in $Y$ used to measure 2-cocycles can be
perturbed to lie in $X$.  Furthermore, if a 2-cycle bounds a 3-chain in $Y$,
the 3-chain can be perturbed to lie in $X$ as well.

Finally in the general case, let $L = \{K_1,K_2,\ldots,K_k\}$ be a
$k$-component link in $M$.  Given an arbitrary propagator $\alpha$ on $\hP(M)$,
we choose 1-cocycles $\beta_1,\ldots,\beta_k$ and construct propagators
$\alpha_i$ and $\alpha_{i,j}$ as above.  If $I$ has at least three elements,
then the dismemberment $\hP(B_I)$ of $\hP(M_I)$ consists entirely of shared
pieces.  We define
$$\alpha_I = \alpha + \sum_{i \in I} \beta_i \tensor \beta_i$$
on $\hP(M\sm L)$.  We extend $\alpha_I$ to each of the other shared pieces by
reusing either $\alpha_i$ or $\alpha_{i,j}$. The conclusion is the following
technical lemma:

\begin{lemma} Let $L \subset M$ be a link with $k$ components. For each
component $K_i$ let $\beta_i$ be a Seifert cocycle. Then for each
$I \subseteq [k]$, there is
cohomological propagator $\alpha_I$ on the manifold $M_I$ such that
$$\alpha_I = \alpha + \sum_{i \in I} \beta_i \tensor \beta_i$$
on $\hP(M \sm L)$, and otherwise $\alpha_I$ and $\alpha_{I'}$ agree on each
component shared by the dismemberments $\hP(B_I)$ 
and $\hP(B_{I'})$ of $\hP(M_I)$ and $\hP(M_{I'})$ along
$\d L$. \label{l:astechnical}
\end{lemma}

Finally let $M_\dis$ be the union of the dismemberments $B_I$ of all $M_I$.
Likewise
let $P(M_\dis)$, $C_{\Gamma,\infty}(M_\dis)$, and $c_\Gamma(E_\dis)$
be the union, respectively, of all dismemberments of topological propagators,
configuration spaces, and bundles associated to each $M_I$.
We extend each $\alpha_I$ by 0 to define it on all of $P(M_\dis)$.

For Torelli surgery we will use a dual construction called \emph{bubble wrap}
in which we glue configuration spaces together instead of dismembering them.
More precisely, if $T \subset M$ is a multi-handlebody with $k$ components
$H_1,\ldots,H_k$, and if $H'_1,\ldots,H'_k$, the bubble-wrap model $M_\bub$ of
the pair $(M,T)$ is given by gluing in both $H_i$ and $H'_i$ to $M\setminus T$
for each $i$. The topological propagator $\hP(M_\bub)$ and the configuration
space $\hC_n(M_\bub)$ are likewise formed from $M_I$ by gluing together
$\hP(M_I)$ and $\hC_n(M_I)$,  ranging over all $I \subseteq [k]$,
wherever these spaces agree. In the Torelli analogue of
Lemma~\ref{l:difference}, $\alpha_1 - \alpha$ is null-homologous where it is
defined.  Consequently the above reasoning allows us to choose the $\alpha_I$
to agree on their common domains, which leads to the following conclusion.

\begin{lemma} Let $T \subset M$ be a multi-handlebody with $k$ components. 
Then the cohomological propagators $\alpha_I$, ranging
over all $I \subseteq [k]$, form a cocycle
$\alpha \in \hP(M_\bub)$.
\label{l:ttechnical}
\end{lemma}

Since in the bubble wrap model there is only one cocycle, we will instead add
and subtract cycles.  For this purpose, given a weight system $w$, we define
$\mu_{w,I}$ as a cycle on $C_{\Gamma,\infty}(M_\bub)$ by extending $\mu_w$,
which exists on $C_{\Gamma,\infty}(M_I)$, by 0.  Dually, all $\mu_{w,I}$ exist
as chains on their common domains on $C_{\Gamma,\infty}(M_\dis)$ and on
$c_\Gamma(E_\dis)$, although they are no longer cycles because of dismemberment
and because we have suppressed gluing.  They form a chain $\mu_w$.

\subsection{The invariants are finite type}
\label{s:fintype}

\subsubsection{Torelli surgery}

We first discuss the Torelli case since it is a bit simpler than the
algebraically split case.  In light of Lemma~\ref{l:propcohom}, we cannot take
a tensor power of a cohomological propagator $\alpha \in H^2(\hP(M))$ at the
cochain level; instead we use a tensor product
$$\alpha_1 \tensor \alpha_2 \tensor \ldots \tensor \alpha_{3n}.$$
Nonetheless the arguments of Section~\ref{s:gory} apply to each $\alpha_i$
separately.  For brevity we let $\gamma$ be its pull-back under $\Phi^*$ to 
$\hC_{\Gamma,\infty}(M_\bub)$.

The constructions of Section~\ref{s:gory} leave us with a cocycle $\alpha$ on
$\hP(M_\bub)$ as well as a family of cycles $\mu_{w,I}$ on
$\hC_{\Gamma,\infty}(M_\bub)$, and we wish to compute the alternating sum of
pairings
$$I^{(k)}_w(M,T) = \sum_{I \subseteq [k]} (-1)^{|I|}
\langle \mu_{w,I}, \gamma \rangle.$$
Observe that the cycles  $\mu_I$ form a parallelepiped in the vector
space of all cycles on $\hC_{\Gamma,\infty}(M_\bub)$.
In other words, there is a cycle-valued,
affine-linear functional $\mu_w(t)$, where $t \in \R^k$ is a vector
of parameters, such that
$$\mu_{w,I} = \mu_w(t_I),$$
where $(t_I)_i$ is 1 for $i \in I$ and 0 for $i \not\in I$.
Let
$$I_w(t) = \langle \mu_w(t), \gamma \rangle.$$
Then $I^{(k)}_w(M,T)$ is a finite difference
\eq{e:findiff}{I^{(k)}_w(M,T) = \left.\del_1 \del_2 \ldots \del_k
    I_w(t)\right|_0,}
where by definition
$$\del_i I = I(t_i) - I(t_i+1).$$
Also let
$$\nu = \left.\del_1 \del_2 \ldots \del_k \mu_w(t) \right|_0$$
be the finite difference as the cycle level; then
\eq{e:newpairing}{I^{(k)}_w(M,T) = \langle \nu_w, \gamma \rangle}

Equation~\ref{e:findiff} passes from formal finite differences of 3-manifold
invariants to traditional finite differences of polynomials. It follows that
$I^{(k)}_w(M)$ vanishes when $k > 2n$, because $I_w(t)$ is a polynomial of
degree $2n$ in $t$. Indeed, the cycle-valued finite difference $\nu$
vanishes identically when $k > 2n$.

A more precise calculation gives us the borderline finite
difference $I^{(2n)}_w(M,T)$.  Observe that
if a configuration $f:\Gamma \to M_\bub$ is disjoint
from a bubble
$$B_i = H_i \cup H'_i$$
of the Torelli surgery then at this point $\mu_w(t)$ is independent of $t_i$;
consequently $\nu$ vanishes here.  Since there are as many bubbles as vertices,
$\Gamma$ must have exactly one vertex in each bubble in the non-vanishing part
of the pairing.  Moreover the bubbles are 3-manifolds; on their product, the
cycle $\nu_w$ is just the fundamental homology class times the weight
$w(\Gamma)$. So we may write the pairing \eqref{e:newpairing} as
\eq{e:bubpairing}{\sum_\Gamma w(\Gamma) \sum_{f:\Gamma \to [2n]}  \langle 
[B_1 \times B_2 \times \ldots \times B_{2n}],\gamma \rangle.}
Given that in this sum each edge of $\Gamma$ connects two distinct bubble $B_i$
and $B_j$, the corresponding factor of the cohomological propagator $\gamma$
measures the linking between 1-cycles in the handlebody $H_i$ (or $H'_i$) and
1-cycles in the handlebody $H_j$ (or $H'_j$).   This linking is the same before
and after Torelli surgery, and the inclusion $H_i \subset B_i$ is an
isomorphism in first homology.

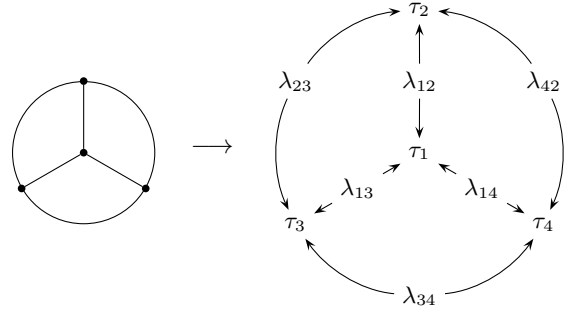
\begin{fullfigure}{f:arrow}{Replacing a Jacobi diagram by a contracted tensor}
\pspicture[.5](-2,-2)(2.25,2)
\cnode*(0,0){.08}{a}\cnode*(1.5;90){.08}{b}
\cnode*(1.5;210){.08}{c}\cnode*(1.5;330){.08}{d}
\ncline{a}{b}\ncline{a}{c}\ncline{a}{d}\pscircle(0,0){1.5}
\endpspicture
$\longrightarrow$
\pspicture[.5](-3.75,-3)(3.25,3)
\rput(0,0){\rnode{a}{\small$\tau_1$}}
\rput(3;90){\rnode{b}{\small$\tau_2$}}
\rput(3;210){\rnode{c}{\small$\tau_3$}}
\rput(3;330){\rnode{d}{\small$\tau_4$}}
\rput(1.5;90){\rnode{ab}{\small$\lambda_{12}$}}
\rput(1.5;210){\rnode{ac}{\small$\lambda_{13}$}}
\rput(1.5;330){\rnode{ad}{\small$\lambda_{14}$}}
\rput(3;150){\rnode{bc}{\small$\lambda_{23}$}}
\rput(3;270){\rnode{cd}{\small$\lambda_{34}$}}
\rput(3;30){\rnode{bd}{\small$\lambda_{42}$}}
\psset{nodesep=3pt,arcangle=-25}
\ncline{<-}{a}{ab}\ncline{<-}{a}{ac}\ncline{<-}{a}{ad}
\ncline{<-}{b}{ab}\ncline{<-}{c}{ac}\ncline{<-}{d}{ad}
\ncarc{<-}{b}{bc}\ncarc{->}{bc}{c}
\ncarc{<-}{c}{cd}\ncarc{->}{cd}{d}
\ncarc{<-}{d}{bd}\ncarc{->}{bd}{b}
\endpspicture
\end{fullfigure}

In conclusion the pairing \eqref{e:bubpairing} becomes a contraction of
tensors:  A vertex in the bubble $B_i$ is replaced by the trilinear form
$$\tau_i:H^1(B_i)^{\times 3} \to \Q$$
given by the triple cup product, an edge connecting $B_i$ to $B_j$ is replaced
by the pairing
$$\lambda_{i,j}:H_1(B_i) \times H_1(B_j) \to \Q$$
given by linking in
any $M_I$, and when an edge is incident to a vertex, the tensors are
contracted.  \fig{f:arrow} gives an example of such a replacement
using arrow notation for tensor
contractions~\cite{Kuperberg:involutory}. These tensor expressions are
summed over Jacobi diagrams $\Gamma$ with vertices decorated by
bubbles.  Finally there is a factor of $2^{3n}(3n)!$ arising from
orderings and orientations of the edges of $\Gamma$, which are now
vestigial.  This leads to the desired value for $I^{(2n)}(M,T)$
(implicit in work of Garoufalidis and Levine~\cite{GL:blinks}).

\subsubsection{Algebraically split surgery}
\label{s:as}

In the algebraically split case, there is one chain $\mu_w$ on all of
$C_{\Gamma,\infty}(M_\dis)$ and on $c_{\Gamma}(E_\dis)$, but there are $2^k$
cocycles $\alpha_I$.  These also form a parallelepiped in the space of all
cocycles on $\hP(M_\dis)$, which is also encoded by an affine-linear function
$\alpha(t)$, with $t \in \R^k$, such that 
$$\alpha_I = \alpha(t_I)$$
for all $I \subseteq [k]$.  By Lemma~\eqref{l:astechnical}, the
function $\alpha(t)$ has the
explicit form
\eq{e:linear}{\alpha(t) = \alpha + \sum_i t_i \beta_i \tensor \beta_i}
on the link complement $M\setminus L$; slightly more generally, the formula
also shows the dependence of $\alpha(t)$ on $t_i$ everywhere outside of the
component $K_i$. We correspondingly let
$$\gamma(t) = \Phi^*(\alpha_1(t) \tensor \alpha_2(t) \tensor \ldots \tensor \alpha_k(t)$$
in keeping with Lemma~\ref{l:propcohom}, and we define
$$\kappa = \left. \del_1 \del_2 \ldots \del_k \gamma(t)\right|_0.$$
We would like to compute
\eq{e:asnewpairing}{I^{(k)}(M,L) = \sum_\Gamma \langle \mu_w, \kappa \rangle.}
Two properties of this finite difference can be argued relatively
easily.  If $k>3n$, then $\kappa$ vanishes identically, because
$\gamma(t)$ is a polynomial of degree $3n$ in $t$.  If $k>2$,
then $\kappa$ vanishes on $c_{\Gamma}(E_\dis)$, because on
each component of $E_\dis$, $\gamma(t)$ is either proportional
to a single $t_i$ (if the component bounds the knot $K_i$)
or it is constant (if the component is shared for all surgeries).

As with Torelli surgery, the marginal case $k=3n$ simplifies because  $\kappa$
is non-zero only when $\Gamma$ is distributed among all components of the
surgery. The following lemma expresses this principle of resource exhaustion.

\begin{lemma}\label{l:resource} If $k=3n$, then 
$\gamma(t)$ is independent of some $t_i$ at a configuration 
$f:\Gamma \to M_\dis$ unless each point in $f(\Gamma)$ lies at
a triple intersection of Seifert surfaces of the link $L$ in $M$.
\end{lemma}
\begin{proof} Say that a vertex of $\Gamma$ provides a dollar to the component
$K_i$ if it lies in the knot $K_i$ in $L$, and that it provides 50 cents if it
lies in the Seifert surface $S_i$.  By equation \eqref{e:linear}, each component $K_i$,
of which there are $3n$, needs a dollar in order for $\kappa(t)$ to depend on
$t_i$ at the configuration $f$. Each vertex, of which there are $2n$, can
provide at most \$1.50, and only by lying at the intersection of three Seifert
surfaces.  The components need $3n$ dollars, which is the most that the
vertices can provide.  Therefore the vertices lie on the Seifert surfaces.
\end{proof}

Having established that the finite difference $\kappa$ is supported
in the link complement $M\setminus L$, we can compute $I^{(3n)}(M,L)$
using the relative cohomology ring $H^*(M,L)$.   Equation~\eqref{e:linear}
implies that
$$\kappa = \Phi^*(\prod_i \beta_i \tensor \beta_i).$$
This cocycle blows down from the configuration space $C_{\Gamma,\infty}(M\setminus L)$
to the Cartesian product $(M \setminus L)^{\Gamma}$.
After blowing down, the chain $\mu_w$ is now proportional to 
the fundamental class:
$$\mu_w = w(\Gamma)[(M \setminus L)^{\Gamma}].$$
The upshot is that the pairing~\eqref{e:asnewpairing} evaluates to another
numerical formula with the geometry of $\Gamma$:  the total weight of all
diagrams $\Gamma$ decorated with a bijection with the link components. Here the
weight of any single diagram is the product of the weights of its vertices.  If
a vertex has incoming edges $i$, $j$, and $k$, its weight is the triple linking
number of the knots $K_i$, $K_j$, and $K_k$.  This is again the desired
answer~\cite{GL:blinks}.

\begin{remark}
Blowing down from the configuration space to the Cartesian product is one
solution to a geometric difficulty in the computation of $I^{(k)}(M,L)$:  Two
vertices of $\Gamma$ might want to lie at the same triple intersection of
Seifert surfaces in $M\setminus L$,  but it is then difficult to see the
behavior of the propagator between them.  In differential terms, the operation
of blowing down says that the diagonal singularities of the propagators cancel
when we take suitable finite differences.  Another approach is to choose two
Seifert surfaces $S_i$ and $S'_i$ for each link component $K_i$, so that
Lemma~\eqref{l:astechnical} becomes
$$\alpha_I = \alpha + \sum_i \beta_i \tensor \beta'_i.$$
If all of the Seifert surfaces are in general position, then the triple points
on $S_i$ and on $S'_i$ will be disjoint, and the computation of $I^{(3n)}(M,L)$
reduces to counting transverse intersections of manifolds far away from the
blowup loci.
\end{remark}

For rational homology spheres there is an interesting generalization of
algebraically split surgery: the framing of each link component $K_i$ can be a
non-zero rational number $p_i/q_i$. In this case Lemma~\ref{l:astechnical}
becomes
$$\alpha_I = \alpha + \sum_i \frac{q_i}{p_i} \beta_i \tensor \beta_i.$$
It follows that the marginal finite difference $I^{(3n)}(M,L)$
is multilinear in the reciprocals of the framings.

\subsection{An unframed invariant}
\label{s:unframed}

The proof of Theorem~\ref{th:correction} rests on three constructions.

First, let $W$ be a closed homology 4-manifold with a 3-plane bundle $E$.
Following Section~\ref{s:bordism}, the sphere bundle $c_e(E)$ has a canonical
cohomology class $\alpha$ which is antisymmetric with respect to fiberwise
inversion, and there is a bundle $c_n(E)$ of
total configuration spaces of the fibers.  As usual, the pull-back
$$\Phi^*(\alpha^{\tensor 3n}) \in H^{6n}(c_n(E))$$
maps to an element in the Jacobi diagram space $V_n$,
yielding a universal invariant $I_n(E)$.

The class $\alpha$ is not only canonical, but \emph{functorial} with respect to
pull-backs of bundles.  The rest of the construction is fiberwise and
therefore also functorial.  On the other hand, since $E$ is a real 3-plane
bundle, its only rational characteristic number is its Pontryagin number
$p_1(E)$ \cite{MS:characteristic}.  Consequently
$$I_n(E) = r_n p_1(E)$$
for some universal vector $r_n \in V_n$.

Second, if $F$ is an oriented 4-plane bundle over some space, it has two
associated 3-plane bundles $F_\pm = \Lambda_2^\pm(F)$ whose fibers are the
spaces of self-dual and anti-self-dual antisymmetric 2-tensors. If $W$ is an
orientable Riemannian 4-manifold with boundary $M$, the bundles $T_\pm W$ both
canonically restrict to $TM$.  Also $W$ has a modified tangent bundle $T'W$
that extends $T'M$, and correspondingly $T'_\pm W$ extend $T'M$.

If $F$ is any oriented 4-plane bundle over it, then the average of
the Pontryagin numbers $p_1(F_\pm)$ is the
Pontryagin number $p_1(F)$.
If $F = TW$, then the Hirzebruch signature theorem
says that the Pontryagin number is thrice the signature $\sigma(W)$, defined as
$a-b$ if the intersection form of $W$ has signature $(a,b)$ \cite[Th.
19.4]{MS:characteristic}.  Algebraically,
\eq{e:Hirz}{p_1(TW) = \frac12(p_1(T_+W) + p_1(T_-W)) = 3\sigma(W).}
If $W = W_1 \cup W_2$ is the union of two 4-manifolds which share boundary $M$
with a marked point $\infty$, then it has a modified tangent bundle $T'W$
extending $T'W_1$ and $T'W_2$.  The Euler number of $T'W$ differs by 2 from
that of $TW$, but the Pontryagin number is the same, so $T'W$ satisfies
equation~\eqref{e:Hirz} as well. (Since $TW$ and $T'W$ differ only in the
neighborhood of $\infty$ and only in a canonical way, this fact can be verified
with a single example: $T'S^4$ is trivial while $\chi(TS^4)=2$ and $p_1(TS^4) =
0$.)

Third, if $M$ is a homology 3-sphere decorated with a bundle bordism $E$ over a
homology 4-manifold $W$.  Let $W_1$ be a smooth 4-manifold with boundary $M$
and signature 0, and let $E_\pm$ be the bundles formed by extending $E$ by
$T'_\pm W_1$.  Then we define
$$\delta_n(M) = \frac{r_n}2(p_1(E^+) + p_1(E^-)).$$
By equation~\ref{e:Hirz}, this quantity does not depend on $W_1$.  (If we
replace $W_1$ by $W_2$, their union has signature 0 because $M$ is a rational
homology sphere; consequently the Pontryagin number, which determines the
change in $\delta_n(M)$, is 0 as well.)  Also the difference
$$\tI_n(M) = I_n(M) - \delta_n(M)$$
is independent of $E$ by the definition of $r_n$.

It remains to show that $\delta_n(M)$ is finite type of degree 1.  The
argument is clearer if we restrict to certain specific bundle bordisms on $M$
and its relatives obtained by surgery.  Namely we choose a 4-manifold $W$ with
boundary $M$ and we decorate $M$ with the formal average of the bundles
$\Lambda_2^\pm(T'W)$.  
In this case the framing
correction is given by
$$\delta_n(M) = 3 r_n \sigma(W).$$
If we perform surgery on a knot $K \in M$ or a Torelli surgery on a handlebody
$H \in M$, we extend $W$ arbitrarily.   In this case the intersection form of
$W$ changes by taking direct sums with matrices that depend only on the
surgery.  Since the signature of a form is linear under direct sums, it is
finite type of degree 1, as desired.  The argument that it is finite type for
 general decorations of $M$ is similar.

To conclude this section, 
we compute the first framing correction coefficient
$r_1$.  The invariant $I_1(M)$ lies in the 1-dimensional vector
space $V_1$ generated by a theta graph; we choose a basis
such that
$$I_1(M) = \langle\alpha^{\cup 3}/6,[\hC_2(M)]\rangle.$$
The simplest twisted $S^2$-bundle
on $S^4$ has Pontryagin number $p_1 = 4$ and total space $\C P^3$; one
model of it is the sequence
$$S^7 \to \C P^3 \to \H P^1 = S^4$$
given by quotienting $S^7$ by complex and quaternionic multiplication.
By the ring structure of $H^*(\C P^3)$ and the fact that
$\alpha$ generates it,
$$\langle\alpha^{\cup 3}/6,[\C P^3]\rangle = 1/6.$$
Thus $r_1 = 1/24$.


\providecommand{\bysame}{\leavevmode\hbox to3em{\hrulefill}\thinspace}

\end{document}